\documentclass[a4paper,10pt]{article}
\usepackage{amsmath,amsthm}
\usepackage{amsfonts}
\usepackage{mathabx}
\usepackage{amssymb,dsfont,graphicx}
\usepackage{caption, subcaption, color,commath} 
\usepackage{cases, xcolor}
\usepackage[left=2.5cm,right=2.5cm,top=3cm,bottom=3cm]{geometry}
\usepackage{hyperref}
\hypersetup{
    colorlinks = true,
    citecolor = blue,
}

\usepackage{todonotes}

\allowdisplaybreaks[1]
\newtheorem{theorem}{Theorem}
\newtheorem{proposition}{Proposition}
\newtheorem{corollary}{Corollary}

\newtheorem{lemma}{Lemma}

\newtheorem{remark}{Remark}

\newcommand{\etal}{\textit{et al.\,}}

\newcommand{\E}{\mathbb{E}}
\newcommand{\id}{\mathds{1}}
\newcommand{\Prob}{\mathbb{P}}

\newcommand{\R}{\mathbb{R}}
\newcommand{\LL}{\mathcal{L}}
\newcommand{\Hn}{\mathcal{H}_n}
\newcommand{\usup}[1]{\,\underset{#1}{\sup}\,}
\newcommand{\umax}[1]{\,\underset{#1}{\max}\,}
\newcommand{\uinf}[1]{\,\underset{#1}{\inf}\,}

\newcommand{\sumi}{\sum_{i=1}^n\,}
\newcommand{\norme}[2]{\left\|#1\right\|_{#2}}
\newcommand{\normsup}[2]{\left\|#1\right\|_{\infty,#2}}
\newcommand{\normsupglobal}[1]{\left\|#1\right\|_{\infty}}
\DeclareMathOperator{\argmin}{argmin\,}

\DeclareMathOperator{\MSE}{MSE}

\DeclareMathOperator{\card}{card}

\newcommand{\red}[1]{\textcolor{red}{#1}}

\author{Ga\"elle Chagny\footnote{LMRS, UMR CNRS 6085, Universit\'e de Rouen Normandie, \texttt{gaelle.chagny@univ-rouen.fr}}, Antoine Channarond\footnote{LMRS, UMR CNRS 6085, Universit\'e de Rouen Normandie, \texttt{antoine.channarond@univ-rouen.fr}}, Van H\`a Hoang\footnote{LMRS, UMR CNRS 6085, Universit\'e de Rouen Normandie, \texttt{van-ha.hoang@univ-rouen.fr}}, Angelina Roche\footnote{Universit\'e Paris-Dauphine, CNRS, UMR 7534, CEREMADE, 75016 Paris, France,  \texttt{roche@ceremade.dauphine.fr}}}

\title{Adaptive nonparametric estimation of a component density in a two-class mixture model}
\date{\today}
\begin{document}
\maketitle

\begin{abstract}
A two-class mixture model, where the density of one of the components is known, is considered. We address the issue of the nonparametric adaptive estimation of the unknown probability density of the second component. We propose a randomly weighted kernel estimator with a fully data-driven bandwidth selection method, in the spirit of the Goldenshluger and Lepski method. An oracle-type inequality for the pointwise quadratic risk is derived as well as convergence rates over H\"older smoothness classes.  The theoretical results are illustrated by numerical simulations.
\end{abstract}

\section{Introduction}\label{sec:intro}
The following mixture model with two components:
\begin{equation}\label{model}
g(x) = \theta + (1-\theta)f(x), \quad \forall x\in [0,1],
\end{equation}
where the mixing proportion $\theta \in (0,1)$ and the probability density function $f$ on $[0,1]$ are unknown, is considered in this article. It is assumed that $n$ independent and identically distributed (\textit{i.i.d.} in the sequel) random variables $X_1,\dots,X_n$ drawn from density $g$ are observed. The main goal is to construct an adaptive estimator of the nonparametric component $f$ and to provide non-asymptotic upper bounds of the pointwise risk : the resulting estimator should automatically adapt to the unknown smoothness of the target function. The challenge arises from the fact that there is no direct observation coming from $f$. As an intermediate step, the estimation of the parametric component $\theta$ is addressed as well.

Model \eqref{model} appears in some statistical settings: robust estimation and multiple testing among others. The one chosen in the present article, as described above, comes from the multiple testing framework, where a large number $n$ of independent hypotheses tests are performed simultaneously. $p$-values $X_1,\dots,X_n$ generated by these tests can be modeled by \eqref{model}. Indeed these are uniformly distributed on $[0,1]$ under null hypotheses while their distribution under alternative hypotheses, corresponding to $f$, is unknown. The unknown parameter $\theta$ is the asymptotic proportion of true null hypotheses. It can be needed to estimate $f$, especially to evaluate and control different types of expected errors of the testing procedure, which is a major issue in this context. See for instance Genovese and Wassermann \cite{genovese2002operating}, Storey \cite{storey2002}, Langaas \etal \cite{langaas2005}, Robin \etal \cite{Robin07}, Strimmer \cite{strimmer2008unified}, Nguyen and Matias \cite{nguyen_matias_2014}, and more fundamentally, Benjamini \etal \cite{Benjamini1995} and Efron \etal \cite{efron2001empirical}.


In the setting of robust estimation, different from the multiple testing one, model \eqref{model} can be thought of as a contamination model, where the unknown distribution of interest $f$ is contaminated by the uniform distribution on $[0,1]$, with the proportion $\theta$. This is a very specific case of the Huber contamination model \cite{huber1965robust}. The statistical task considered consists in robustly estimating $f$ from contaminated observations $X_1,\dots,X_n$. But unlike our setting, the contamination distribution is not necessarily known while the contamination proportion $\theta$ is assumed to be known, and the theoretical investigations aim at providing minimax rates as functions of both $n$ and $\theta$. See for instance the preprint of Liu and Gao \cite{liu2017density}, which addresses pointwise estimation in this framework.

Back to the setting of multiple testing, the estimation of $f$ in model \eqref{model} has been addressed in several works. Langaas \etal \cite{langaas2005} proposed a Grenander density estimator for $f$, based on a nonparametric maximum likelihood approach, under the assumption that $f$ belongs to the set of decreasing densities on $[0,1]$. Following a similar approach, Strimmer \cite{strimmer2008unified} also proposed a modified Grenander strategy to estimate $f$. However, the two aforementioned papers do not investigate theoretical features of the proposed estimators. Robin \etal \cite{Robin07} and Nguyen and Matias  \cite{nguyen_matias_2014} proposed a randomly weighted kernel estimator of $f$, where the weights are estimators of the posterior probabilities of the mixture model, that is, the probabilities of each individual $i$ being in the nonparametric component given the observation $X_i$. \cite{Robin07} proposes an EM-like algorithm, and proves the convergence to an unique solution of the iterative procedure, but they do not provide any asymptotic property of the estimator. Note that their model $g(x)=\theta \phi(x)+(1-\theta)f(x)$, where $\phi$ is a known density, is slightly more general, but our procedure is also suitable for this model under some assumptions on $\phi$. Besides, \cite{nguyen_matias_2014} achieves a nonparametric rate of convergence $n^{-2\beta/(2\beta+1)}$ for their estimator, where $\beta$ is the smoothness of the unknown density $f$. However, their estimation procedure is not adaptive since the choice of their optimal bandwidth still depends on $\beta$.

In the present work, a complete inference strategy for both $f$ and $\theta$ is proposed. For the nonparametric component $f$, a new randomly weighted kernel estimator is provided with a data-driven bandwidth selection rule. Theoretical results on the whole estimation procedure, especially adaptivity of the selection rule to unknown smoothness of $f$, are proved under a given identifiability class of the model, which is an original contribution in this framework. Major results derived in this paper are the oracle-type inequality in Theorem \ref{th:oracle-inequality}, and the rates of convergence over H\"older classes, which are adapted to the control of pointwise risk of kernel estimators, in Corollary \ref{cor:rate-cv-fhat}. 

Unlike the usual approach in mixture models, the weights of the proposed estimator are not estimates of the posterior probabilities. The proposed alternative principle is simple and consists in using weights based on a density change, from the target distribution $f$, which is not directly reachable, to the distribution of observed variables $g$. A function $w$ is thus derived such that $f(x) = w(\theta, g(x))g(x)$, for all $\theta,x\in [0,1]$. This type of link between one of the conditional distribution given hidden variables, $f$, to the distribution of observed variables $g$, is quite remarkable in the framework of mixture models. It is a key idea of our approach, since it implies a crucial equation for controlling the bias term of the risk, see Subsection \ref{sec:kernel-estimator} for more details. This is necessary to investigate adaptivity using the Goldenshluger and Lespki (GL) approach \cite{GL11}, which is known in other various contexts, see for instance, Comte \etal \cite{ComteSamson13}, Comte and Lacour \cite{ComteLacour}, Doumic \etal \cite{DHRR}, Reynaud-Bouret \etal \cite{R-B14} who apply GL method in kernel density estimation, and Bertin \etal \cite{BLR2016}, Chagny \cite{Chagny13}, Chichignoud \etal \cite{Chichignoud2017} or Comte and Rebafka \cite{ComteRebafka2016}.

Thus oracle weights are defined by $w(\theta,g(X_i))$, $i=1,\dots,n$, but $g$ and $\theta$ are unknown. These oracle weights are estimated by plug-in, using preliminary estimators of $g$ and $\theta$, based on an additional sample $X_{n+1},\dots,X_{2n}$. Some assumptions on these estimators are needed to prove the results on the estimator of $f$; this paper also provides estimators of $g$ and $\theta$ which satisfy these assumptions. Note that procedures of \cite{nguyen_matias_2014} and \cite{Robin07} actually require preliminary estimates of $g$ and $\theta$ as well, but they do not deal with additional uncertainty caused by the multiple use of the same observations in the estimates of $\theta$, $g$ and $f$. 

Identifiability issues are reviewed in Section 1.1 in Nguyen and Matias \cite{Ng-Matias-SJS}. In the present work, $f$ is assumed to be vanishing at a neighbourhood of $1$ to ensure identifiability. Under this assumption, $\theta$ can be recovered as the infimum of $g$. Moreover, as shown above by the equation linking $f$ to $g$ and $\theta$, $f$ is actually uniquely determined by giving $g$ and $\theta$, even though the latter is not the infimum of $g$. Note that the theoretical results on the estimator of the nonparametric component $f$ do not depend on the chosen identifiability class, and can be transposed to other cases. For that reason, the discussion on identifiability is postponed to Section \ref{sec:estimation-theta}, after results on the estimator of $f$.



The paper is organized as follows. Our randomly weighted estimator of $f$ is constructed in Section \ref{sec:kernel-estimator}. Assumptions on $f$ and on preliminary estimators of $g$ and $\theta$ required for proving the theoretical results are in this section too. In Section \ref{sec:estimation}, a bias-variance decomposition for the pointwise risk of the estimator of $f$ is given as well as the convergence rate of the kernel estimator with a fixed bandwidth. In Section \ref{sec:adaptation}, an oracle inequality is given, which justifies our adaptive estimation procedure. Construction of the preliminary estimators of $g$ and $\theta$ are to be found in Section \ref{sec:estimation-g-theta}. Numerical results illustrate the theoretical results in Section \ref{sec:simus}. Proofs of theorems, propositions and technical lemmas are postponed to Section \ref{sec:proofs}.
 
\section{Collection of kernel estimators for the target density}\label{sec:estimation}

In this section, a family of kernel estimators for the density function $f$ based on a sample $(X_i)_{i=1,\ldots,n}$ of i.i.d. variables with distribution $g$ is defined. It is assumed that preliminary estimators of both the mixing proportion $\theta$ and the mixture density $g$ are available, and respectively denoted by $\tilde{\theta}_n$ and $\hat g$. They are defined from an additional sample $(X_i)_{i=n+1,\ldots,2n}$ of independent variables also drawn from $g$ but independent of the first sample $(X_i)_{i=1,\ldots,n}$. Definitions, results and results on these preliminary estimates are the subject of Section \ref{sec:estimation-g-theta}.

\subsection{Construction of the estimators}\label{sec:kernel-estimator}

To define estimators for $f$, the challenge is that observations $X_1,\ldots, X_n$ are not drawn from $f$ but from the mixture density $g$. Hence the density $f$ cannot be estimated directly by a classical kernel density estimator. Thus we will build weighted kernel estimates  This idea has been used in other contexts, see for example \cite{ComteRebafka2016}. The starting point is the following lemma whose proof is straightforward.

\begin{lemma}\label{lem:transition-property}
Let $X$ be a random variable from the mixture density $g$ defined by \eqref{model} and $Y$ be an (unobservable) random variable from the component density $f$. Then for any measurable bounded function $\varphi$:
\begin{equation}\label{eq:transition-property}
\E \big[\varphi(Y)\big] = \E \big[w(\theta,g(X))\varphi(X)\big],
\end{equation}
where
\[
w(\theta,g(x)) := \frac 1{1-\theta}\left(1 - \frac{\theta}{g(x)}\right), \,\, x\in [0,1].
\]
\end{lemma}

This result will be used as follows.  Let $K:\R \to \R$ be a  kernel function, that is an integrable function such that $\int_{\mathbb R} K(x)dx = 1$ and $\int_{\mathbb R}  K^2(x)dx < +\infty$. For any $h>0$, let $K_h(\cdot) = K(\cdot/h)/h$. Then the choice $\varphi(\cdot) = K_h(x - \cdot)$ in Lemma \ref{lem:transition-property} gives:
\[
\E \big[K_h(x-Y) \big] = \E \big[w(\theta,g(X))K_h(x-X)\big],
\]

This leads to define the following randomly weighted kernel estimator of $f$:
\begin{equation}\label{eq:estimator-f}
\hat{f}_h(x) = \frac 1n\sumi w(\tilde \theta_n, \hat{g}(X_i))K_h(x-X_i),\,\, x\in [0,1],
\end{equation}
where:
\begin{equation}\label{eq:weights}
w(\tilde \theta_n, \hat{g}(X_i)) = \frac 1{1-\tilde \theta_n}\left(1 - \frac{\tilde \theta_n}{\hat{g}(X_i)}\right),\quad i=1,\ldots, n. 
\end{equation}
Therefore, $\hat{f}_h$ is a randomly weighted kernel estimator of $f$. Note that the total sum of the weights may not equal $1$, in comparison with the estimators proposed in Nguyen and Matias \cite{nguyen_matias_2014} and Robin \etal \cite{Robin07}. The main advantage of such weights, is that, if we replace $\hat g$ and $\tilde{\theta}_n$ by their theoretical unknown counterparts $g$ and $\theta$ in \eqref{eq:estimator-f}, we obtain, $\mathbb E[\hat{f}_h(x)]=K_h\star f(x)$, where $\star$ stands for the convolution product. This relation, classical in nonparametric kernel estimation, is crucial to study the bias term in the risk of the estimator, and hence to reach adaptivity.

\subsection{Risk bounds of the estimator}\label{sec:BV-decomp-fixed-bw}

Here, upper bounds are derived for the pointwise mean-squared error of the estimator $\hat f_h$, defined in \eqref{eq:estimator-f}, with a fixed bandwidth $h>0$.  Our objective is to study the pointwise risk for the estimation of the density $f$ at a point $x_0\in [0,1]$. Throughout the paper, the kernel $K$ is chosen compactly supported on an interval  $[-A,A]$ with $A$ a positive real number, and such that $\sup_{x\in [-A,A]}|K(x)|<\infty$. We denote by $\mathcal{V}_n(x_0)$ the neighbourhood of $x_0$ used in the sequel and defined by
\[
\mathcal{V}_n(x_0) = \left[x_0 - \frac{2A}{\alpha_n}, x_0 + \frac{2A}{\alpha_n} \right],
\]
where $(\alpha_n)_n$ is a positive sequence of numbers larger than 1, only depending  on $n$ such that $\alpha_n \to +\infty$ as $n \to +\infty$, chosen by the user.  For any function $u$ on $\mathbb R$, and any interval $I\subset\R$, let $\normsup{u}{I} = \sup_{t\in I}|u(t)|$. We also denote by  $\gamma = \underset{t\in \mathcal{V}_n(x_0)}{\inf} |g(t)|$. Thanks to \eqref{model}, we have $g(t)\geq \theta>0$ for any $t\in [0,1]$, and thus, $\gamma>0$.

\bigskip

In the sequel, we consider the following assumptions. Note that all assumptions are not simultaneously necessary for the results.

\begin{itemize}
\item[{\bf (A1)}] The density $f$ is uniformly bounded on $\mathcal{V}_n(x_0)$ for some $n$: $\normsup{f}{\mathcal{V}_n(x_0)} < \infty$.\smallskip
\item[{\bf (A2)}] The preliminary estimator $\hat g$ is bounded away from $0$ on $\mathcal{V}_n(x_0)$ a.s. : 
\begin{equation}\label{eq:inf_ghat}
\hat\gamma := \uinf{t\in \mathcal{V}_n(x_0)} |\hat g(t)| >0.
\end{equation}

\item[{\bf (A3)}] The preliminary  estimate $\hat g$ of $g$ satisfies,  for all $\nu >0$ 
\begin{equation}\label{eq:condtion-on-ghat}
\Prob\left(\usup{t\in \mathcal{V}_n(x_0)} \left| \frac{\hat g(t) - g(t)}{\hat g(t)}\right| > \nu \right) \le C_{g,\nu} \exp\left\{ - (\log n)^{3/2} \right\},
\end{equation}
with $C_{g,\nu}$ a constant only depending on $g$ and $\nu$. 

\item[{\bf (A4)}]  The preliminary estimator $\tilde{\theta}_n$ is constructed such that $\tilde{\theta}_n \in [\delta/2, 1-\delta/2]$ a.s., for a fixed $\delta \in (0,1)$.

\item[{\bf (A5)}] For any bandwidth $h>0$, we assume that a.s. 
\[
\alpha_n \le \frac{1}{h} \quad \text{ and } \quad \frac{1}{h} \le \min\left\{\frac{\hat\gamma n}{\log^3(n)} , \frac1n\right\}.
\]
\item[\bf (A6)] $f$ belongs to the H\"older class of smoothness $\beta$ and radius $\LL$ on $[0,1]$, defined by
\[
\Sigma(\beta,\LL) = \left\{\phi: \phi \text{ has } \ell = \lfloor \beta \rfloor \text{ derivatives and } \forall x,y\in [0,1], |\phi^{(\ell)}(x) - \phi^{(\ell)}(y)| < \LL|x - y|^{\beta -\ell}\right\}, 
\]

where $\lfloor x \rfloor$  denotes a smallest integer which is strictly smaller than the real number $x$.
\item[\bf (A7)] $K$ is a kernel of order $\ell $ :  $\int_{\mathbb R} x^j K(x)dx = 0$ for $1\le j \le \ell$ and $\int_{\mathbb R}|x|^{\ell}|K(x)|dx < \infty$.
\end{itemize}

Since $g = \theta + (1-\theta)f$, Assumption {\bf (A1)} implies that $\normsup{g}{\mathcal{V}_n(x_0)} < \infty$. This assumption is needed to control the variance term, among others, of the bias-variance decomposition of the risk. Let us notice that the density $g$ is automatically bounded from below by a positive constant in our model \eqref{model}.  Assumption {\bf (A2)} is required to bound the term $1/\hat{g}(\cdot)$ that appears in the weight $w(\tilde \theta_n, \hat{g}(\cdot))$, see \eqref{eq:weights}.  Assumption {\bf (A3)} means that the preliminary $\hat g$ has to be rather accurate. Assumptions {\bf (A2)} and {\bf (A3)} are also introduced by Bertin \etal \cite{BLR2016} for conditional density estimation purpose : see (3.2) and (3.3) p.946. The methodology used in our proofs is close to their work : the role played by $g$ here corresponds to the role played  by the marginal density of their paper.  They have also shown that an estimator of $g$ satisfying these properties can be built, see Theorem 4, p. 14 of \cite{BLR2013} and some details at Section \ref{sec:estimation-g}. We also build an estimator $\tilde\theta_n$ that satisfies Assumption {\bf (A4)} in Section \ref{sec:estimation-theta}. Assumption {\bf (A5)} deals with the order of magnitude of the bandwidths and is also borrowed from \cite{BLR2016} (see Assumption (CK) p.947). An example of bandwidth collection satisfying Assumption {\bf (A5)} is given in the statement of Corollary~\ref{cor:rate-cv-fhat}. Assumptions {\bf (A6)} and {\bf (A7)} are classical for kernel density estimation, see \cite{Tsybakov04} or \cite{comte_estimation_2015}. The index $\beta$ in Assumption {\bf (A6)} is a measure of the smoothness of the target function. Such assumptions permit to control the bias term of the bias-variance decomposition of the risk, and thus to derive convergence rates.  We will classically choose $\ell=  \lfloor \beta \rfloor$ for Assumption {\bf (A7)} in Corollary \ref{cor:rate-cv-fhat} below.

We first state an upper bound for the pointwise risk of the estimator $\hat{f}_h$. The proof can be found in Section \ref{proof:prop-upper-bound-fhat}.

\begin{proposition}\label{prop:upper-bound-fhat}
Assume that Assumptions {\bf (A1)} to {\bf (A5)} are satisfied. Then, for any $x_0\in [0,1]$ and $\delta \in (0,1)$, the estimator $\hat{f}_h$ defined by \eqref{eq:estimator-f} satisfies
\begin{multline}
\E \left[\big(  \hat f_h(x_0) - f(x_0)\big)^2 \right] \le C^*_1 \left\{ \normsup{K_h\star f - f}{\mathcal{V}_n(x_0)}^2 + \frac{1}{\delta^2\gamma^2 nh}\right\}  \\
+ \frac{C^*_2}{\delta^6} \E \left[\big|\tilde\theta_n - \theta\big|^2 \right] +\frac{C_3^*}{\delta^2\gamma^2}\E  \Big[ \normsup{\hat g - g}{\mathcal{V}_n(x_0)}^2 \Big]   + \frac{C_4^*}{n^2}, \label{eq:upper-bound-fhat}
\end{multline}
where $C^*_{\ell}$, $\ell=1,\ldots,4$ are positive constants such that : $ C^*_1$ depends on $\|K\|_2$ and $\normsup{g}{\mathcal{V}_n(x_0)}$, $C^*_2$ depends on $\normsup{g}{\mathcal{V}_n(x_0)}$ and $\norme{K}{1}$, $C^*_3$ depends on $\norme{K}{1}$ and $C^*_4$ depends on $\normsup{f}{\mathcal{V}_n(x_0)}$, $g$, $\delta$, $\gamma$, and $\normsupglobal{K}$.
\end{proposition}

Proposition \ref{prop:upper-bound-fhat} is a bias-variance decomposition of the risk. The first term in the right-hand-side (\textit{r.h.s.} in the sequel) of \eqref{eq:upper-bound-fhat} is a bias term which decreases when the bandwidth $h$ vanishes whereas the second one corresponds to the variance term and increases when $h$ vanishes. 

There are two additional terms $\E [ \normsup{\hat g - g}{\mathcal{V}_n(x_0)}^2 ]$  and $\E [|\hat \theta_n - \theta|^2 ]$ in the \textit{r.h.s.} of \eqref{eq:upper-bound-fhat}. They are unavoidable since the estimator $\hat{f}_h$ depends on the plug-in estimators $\hat{g}$ and $\tilde{\theta}_n$. However, as proved in Corollary~\ref{cor:rate-cv-fhat}, these two terms does not deteriorate the convergence rate provided that $g$ and $\theta$ are estimated accurately. We define in Section~\ref{sec:estimation-g-theta} such estimators of $g$ and $\theta$. The term $C^*_4/(\delta^2 n^2)$ is a remaining term and is also negligible.

\section{Adaptive pointwise estimation}\label{sec:adaptation}

Let $\mathcal{H}_n$  be a finite family of possible bandwidths $h>0$, whose cardinality is bounded by the sample size $n$. The best estimator in the collection $(\hat f_h)_{h\in\mathcal H_n}$ defined in \eqref{eq:estimator-f} at the point $x_0$ is the one that have the smallest risk, or similarly, the smallest bias-variance decomposition. But since $f$ is unknown, in practice it is impossible to minimize over $\mathcal H_n$ the r.h.s. of inequality \eqref{eq:upper-bound-fhat} in order to select the best estimate. Thus, we propose a data-driven selection, with a rule in the spirit of   Goldenshluger and Lepski (GL in the sequel) \cite{GL11}. The idea is to mimic the bias-variance trade-off for the risk, with empirical counterparts for the unknown quantities. We first estimate the variance term of the trade-off by setting, for any $h\in\Hn$ 
\begin{equation}\label{eq:kappa-variance-term}
V(x_0, h) = \frac{\kappa \norme{K}{1}^2\norme{K}{2}^2\normsup{g}{\mathcal{V}_n(x_0)}}{\hat\gamma^2 nh}\log(n),
\end{equation}
with $\kappa >0$ a tuning parameter. The principle of the GL method is then to estimate the bias term $\normsup{K_h\star f - f}{\mathcal{V}_n(x_0)}^2$ of $\hat f_h(x_0)$ for any $h\in\mathcal H_n$ with 
\begin{equation}\label{eq:def_A}
A(x_0,h) := \umax{h'\in \Hn}\,\left\{ \big( \hat f_{h,h'}(x_0) - \hat f_{h'}(x_0) \big)^2 - V(x_0, h') \right\}_+,
\end{equation}
where, for any $h, h'\in\Hn$,
\[
\hat f_{h,h'}(x_0) = \frac 1n \sumi w(\tilde\theta_n, \hat g(X_i))(K_h\star K_{h'})(x_0 - X_i)  = (K_{h'}\star \hat f_h)(x_0). 
\]
Heuristically, since $\hat{f}_h$ is an estimator of $f$ then $\hat f_{h,h'} = K_{h'}\star \hat  f_h$ can be considered as an estimator of $K_{h'}\star f$. The proof of Theorem \ref{th:oracle-inequality} below in Section \ref{sec:proof_thm} then justifies that $A(x_0,h)$ is a good approximation for the bias term of the pointwise risk. Finally, our estimate at the point $x_0$ is
\begin{equation}\label{eq:adaptive-estimator-f}
\hat f(x_0) := \hat f_{\hat h(x_0)}(x_0),
\end{equation}
where the bandwidth $\hat h(x_0) $ minimizes the empirical bias-variance decomposition :
\[
\hat h(x_0) := \underset{h\in\Hn}{\argmin}\left\{ A(x_0,h) + V(x_0,h) \right\}. 
\]

The constants that appear in the estimated variance $V(x_0,h)$ are known, except $\kappa$, which is a numerical constant calibrated by simulation (see practical tuning in Section \ref{sec:simus}), and except $\normsup{g}{\mathcal{V}_n(x_0)}$, which is replaced by an empirical counterpart in practice (see also Section \ref{sec:simus}). It is also possible to justify the substitution from a theoretical point of view, but it adds cumbersome technicalities. Moreover, the replacement does not change the result of Theorem \ref{th:oracle-inequality} below. We thus refer to Section 3.3 p.1178 in \cite{comte_adaptive_2011} for example, for the details of a similar substitution. The risk of this estimator is controlled in the following result. 

\begin{theorem}\label{th:oracle-inequality}
Assume that Assumptions {\bf (A1)} to  {\bf (A4)} are fulfilled, and that all $h\in\mathcal H_n$ satisfies {\bf (A5)}. Suppose in addition that the sample size $n$ is larger than a constant that only depends on the kernel $K$. For any $\delta \in (0,1)$, the estimator $\hat f(x_0)$ defined in \eqref{eq:adaptive-estimator-f} satisfies 
\begin{eqnarray}
\nonumber\E \left[\big(  \hat f (x_0) - f(x_0)\big)^2 \right] &\leq& C_5^* \underset{h\in \Hn}{\min}\left\{\normsup{K_h\star f - f}{\mathcal{V}_n(x_0)}^2 + \frac{\log(n)}{\delta^2 \gamma^2 nh}\right\}  \\
&&+ \frac{C_6^*}{\delta^6} \usup{\theta \in [\delta, 1- \delta]}\E \left[\big|\tilde\theta_n - \theta\big|^2 \right] + \frac{C_7^*}{\delta^2\gamma^2}\E  \Big[ \normsup{\hat g - g}{\mathcal{V}_n(x_0)}^2 \Big]   + \frac{C_8^*}{n^2}, \label{eq:oracle-inequality}
\end{eqnarray}
where  $C^*_{\ell}$, $\ell=5,\ldots,8$ are positive constants such that  : $C^*_5$ depends on $\normsup{g}{\mathcal{V}_n(x_0)}, \norme{K}{1}$ and $\norme{K}{2}$, $C^*_6$ depends on $\norme{K}{1}$, $C^*_7$ depends on $\normsup{g}{\mathcal{V}_n(x_0)}$ and $\norme{K}{1}$, and $C^*_8$ depends on $\delta$, $\gamma$, $\normsup{f}{\mathcal{V}_n(x_0)}$, $g$, $\norme{K}{2}$ and $\normsupglobal{K}$.
\end{theorem}

Theorem \ref{th:oracle-inequality} is an oracle-type inequality. It holds whatever the sample size, larger than a fixed constant. It shows that the optimal bias variance trade-off is automatically achieved: the selection rule permits to select in a data-driven way  the best estimator in the collection of estimators $(\hat f_h)_{h\in\mathcal{H}_n}$, up to a multiplicative constant $C_5^*$. The last three remainder terms in the \textit{r.h.s.} of \eqref{eq:oracle-inequality} are the same as the ones in Proposition \ref{prop:upper-bound-fhat}, and are unavoidable, as aforementioned. We have an additional logarithmic term in the second term of the \textit{r.h.s.}, compared to the analogous term in \eqref{eq:upper-bound-fhat}. It is classical in adaptive pointwise estimation (see for example \cite{ComteRebafka2016} or \cite{Butucea2000}). In our framework, it does not deteriorate the adaptive convergence rate, see Section \ref{sec:cv_rate} below. To compute this rate, we now have to define estimators for the mixing density $g$ and proportion $\theta$, in such a way that the convergence rate which would be obtained by the minimisation of the first term in the \textit{r.h.s} of \eqref{eq:oracle-inequality} can be preserved. 

\section{Estimation of the mixture density $g$ and the mixing proportion $\theta$}\label{sec:estimation-g-theta}

This section is devoted to the construction of the preliminary  estimators $\hat{g}$ and $\tilde{\theta}_n$,  required to build \eqref{eq:estimator-f}. To define them, we assume that we observe an additional sample $(X_i)_{i={n+1,\ldots,2n}}$ distributed with density function $g$, but independent of the sample $(X_i)_{i=1,\ldots,n}$. We explain how  estimators $\hat g$ and $\tilde \theta_n$ can be defined to satisfy the assumptions described at the beginning of Section \ref{sec:BV-decomp-fixed-bw}, and also how we compute them in practice. The reader should bear in mind that other constructions are possible, but our main objective is the adaptive estimation of the density $f$. Thus, further theoretical studies are beyond the scope of this paper. 

\subsection{Preliminary estimator for the mixture density $g$}\label{sec:estimation-g}

As already noticed, the role played by $g$ to estimate $f$ in our framework finds an analogue in the work of Bertin \textit{et al.} \cite{BLR2016} : the authors propose a conditional density estimation method that involves a preliminary estimator of the marginal density of a couple of real random variables. The assumptions {\bf (A2)} and {\bf (A3)} are borrowed from their paper. From a theoretical point of view, we thus also draw inspiration from them to build $\hat g$.

Since we focus on kernel methods to recover $f$, we also use kernels for the estimation of $g$.  Let $L:\R \to \R$ be a  function such that $\int_{\R} L(x)dx =1$ and $\int_{\R} L^2(x)dx <\infty$. Let $L_b(\cdot) = b^{-1}L(\cdot/b)$, for any $b>0$. The function $L$ is a kernel, but can be chosen differently from the kernel $K$ used to estimate the density $f$. The classical kernel density estimate for $g$ is 
\begin{equation}\label{eq:estim_g}
\hat{g}_{b}(x_0) = \frac{1}{n}\sum_{i=n+1}^{2n} L_b(x_0 - X_i),
\end{equation}
Theorem 4 p.14 of \cite{BLR2013} proves that it is possible to select an adaptive bandwidth $b$ of $\hat{g}_b$ in such a way that  Assumptions {\bf (A2)} and {\bf (A3)} are fulfilled, and that the resulting estimate $\hat g_{\hat b}$ satisfies 
\begin{equation}\label{eq:rate-cv-g}
\E  \Big[ \normsup{\hat g_{\hat b} - g}{\mathcal{V}_n(x_0)}^2 \Big]\leq C\left( \frac{\log n}{n}\right)^{\frac{2\beta}{2\beta + 1}},
\end{equation}
if $g\in\Sigma(\beta,\mathcal L')$, where $C,\mathcal L'>0$ are some constants, and if the kernel $L$ has an order $\ell =  \lfloor \beta \rfloor$. The idea of the result of Theorem 4  in \cite{BLR2013} is to select the bandwidth $\hat b$ with a classical Lepski method, and to apply results from Gin\'e and Nickl \cite{gine_exponential_2009}.  Notice that, in our model, Assumption {\bf (A6)} permits to obtain directly the required smoothness assumption, $g\in\Sigma(\beta,\mathcal L')$. This guarantees that both the assumptions  {\bf (A2)} and {\bf (A3)} on $\hat g$ can be satisfied and that the additional term $\E  [ \normsup{\hat g - g}{\mathcal{V}_n(x_0)}^2]$ can be bounded as required in the statement of Corollary \ref{cor:rate-cv-fhat}.

For the simulation study below now, we start from the kernel estimators $(\hat g_{b})_{b>0}$ defined in \eqref{eq:estim_g} and rather use a procedure in the spirit of the pointwise GL method to automatically select a bandwidth $b$.   First, this choice permits to be coherent with the selection method chosen for the main estimators $(\hat f_{h})_{h\in\mathcal H_n}$, see Section \ref{sec:adaptation}. Then, the construction also provides an accurate estimate of $g$, see for example \cite{comte_estimation_2015}.  Let $\mathcal B$ be a finite family of bandwidths. For any $b, b'\in \mathcal B$, we introduce the auxiliary functions $\hat{g}_{b,b'}(x_0) = n^{-1}\sum_{i=n+1}^{2n} (L_b\star L_{b'})(x_0 - X_i).$ Next, for any $b \in \mathcal B$, we set
\[
A^g(b,x_0) = \underset{b'\in \mathcal B}{\max}\left\{ \left(\hat{g}_{b,b'}(x_0) - \hat{g}_{b'}(x_0) \right)^2  - \Gamma_1(b') \right\}_{+},
\]
where $\Gamma_1(b) = \varepsilon \norme{L}{1}^2\norme{L}{2}^2\normsupglobal{g} \log(n)/(nb)$, with $\varepsilon>0$ a constant to be tuned. Then, the final estimator of $g$ is given by $\hat{g}(x_0) := \hat{g}_{\hat b_g (x_0)}(x_0)$, with $\hat b_g(x_0) := {\argmin}_{b \in \mathcal B}\{ A^g(b ,x_0) + \Gamma_1(b) \}$. The tuning of the constant $\varepsilon$ is presented in Section \ref{sec:simus}.

\subsection{Estimation of the mixing proportion $\theta$}\label{sec:estimation-theta}

A huge variety of methods have been investigated for the estimation of the mixing proportion $\theta$ of model  \eqref{model} : see, for instance, \cite{storey2002},  \cite{langaas2005},  \cite{Robin07}, \cite{celisse2010}, \cite{Ng-Matias-SJS} and references therein.  A common and performant estimator is the one proposed by Storey \cite{storey2002}:  $\theta$ is estimated by $\hat\theta_{\tau,n} = \#\{X_i > \tau;i = n+1,\ldots,2n\}/(n(1-\tau))$ with $\tau$ a threshold to be chosen. The optimal value of $\tau$ is calculated with a boostrap algorithm. However, it seems difficult to obtain theoretical guarantees on $\hat\theta_{\tau,n}$.

For a detailed  discussion about possible identifiability conditions of model \eqref{model}, we refer to Celisse and Robin \cite{celisse2010} or Nguyen and Matias \cite{Ng-Matias-SJS}.  In the sequel we focus on a particular case of model \eqref{model}, which ensures the identifiability of the parameters $(\theta,f)$  (see for example Assumption A in \cite{celisse2010}, or Section 1.1 in \cite{Ng-Matias-SJS}). The density $f$ is assumed to belong to the family 
\begin{equation}\label{eq:identifiability-condition}
\mathcal{F}_\delta = \left\{f:[0,1]\to\R_+,f \text{ is a density such that } f_{|[1-\delta,1]} = 0\right\},
\end{equation}
where $\delta\in (0,1)$. Under this assumption, the main idea to recover $\theta$  is that it is  the lower bound of the density $g$ in model \eqref{model} : $\theta = \inf_{x\in [0,1]}  g(x) = g(1)$. Celisse and Robin \cite{celisse2010} or Nguyen and Matias \cite{Ng-Matias-SJS}  then define a histogram-based estimator $\hat g$ for $g$, and estimate $\theta$ with the lower bound of $\hat g$, or with $\hat g(1)$. The procedure we choose is still based on the same assumption, but, to be consistent with the other estimates, we use kernels to recover $g$ instead of histograms. 

Nevertheless, since it is well-known that kernel density estimation methods suffer from boundary effects, which cause inaccurate estimate of $g(1)$, we cannot directly use the kernel estimates of $g$ defined in \eqref{eq:estim_g}. To deal with this issue, we apply a simple reflection method (see for example Schuster \cite{Schuster85}). From the random sample $X_{n+1},\ldots, X_{2n}$ from density $g$,  we introduce, for $i=1,\ldots,n$,
\begin{equation}\label{eq:Y_i}
Y_i = \begin{cases}
X_{i+n} &\text{ if } \varepsilon_i = 1,\\
2- X_{i+n} &\text{ if } \varepsilon_i = -1,\
\end{cases}
\end{equation}
where $\varepsilon_1, \ldots, \varepsilon_n$ are $n$ \textit{i.i.d.} random variables drawn from Rademacher distribution with parameter $1/2$, and  independent of the $X_i$'s. The random variables $Y_1, \ldots, Y_n$ can be regarded as randomly symmetrized version of the $X_i$'s, with support  $[0,2]$ (see the first point of Lemma \ref{lem:g_sym} below).  Now, suppose that $L$ is a symmetric kernel. For any $b>0$, define
\begin{equation}\label{eq:symmetrized-kernel-theta}
\hat{g}^{sym}_b(x) = \frac{1}{n} \sum_{k=1}^{n}L_b(x - Y_k), \quad x \in [0,2].  
\end{equation}
Instead of evaluating $\hat{g}^{sym}_{b}$ at the single point $x=1$, we evaluate twice the average of all the values of the estimator $\hat{g}^{sym}_{b}$ on the interval $[1-\delta,1+\delta]$, relying on the fact that $\theta = g(x)$, for all  $x \in [1-\delta,1]$ (under the assumption $f\in\mathcal F_{\delta}$), to increase the accuracy of the resulting estimate.  Thus, we set  
\begin{equation}\label{eq:theta-estimator}
\hat\theta_{n,b} = \frac{1}{\delta} \int_{1-\delta}^{1+\delta} \hat{g}^{sym}_{b}(x)dx.
\end{equation}
Finally, for the estimation of $f$, we use a truncated estimator $\tilde{\theta}_n$ defined as
\begin{equation}\label{eq:estimator-tilde-theta}
\tilde\theta_{n,b} := \max\big(\min(\hat\theta_{n,b}, 1- \delta	/2), \delta/2 \big).
\end{equation}
The definition of $\tilde{\theta}_{n,b}$ permits to ensure that $\tilde\theta_{n,b} \in [\delta/2, 1-\delta/2]$ : this is Assumption {\bf (A4)}. This permits to avoid possible difficulties in the estimation of $f$ when $\hat\theta_{n,b}$ is close to zero, see \eqref{eq:estimator-f}. The following lemma establishes some properties of all these estimates. Its proof can be found in Section \ref{sec:proof_lem_g_sym}.

\begin{lemma}\label{lem:g_sym}\quad\vspace*{-0.025in}
\begin{itemize}
\item The random variables $Y_k$, $k\in\{1,\ldots,n\}$, are \textit{i.i.d.}, with density  $$g^{sym}\,:x\; \longmapsto \left\{\begin{array}{l}g(x)/2\mbox{ if }x\in[0,1]\\g(2-x)/2\mbox{ if }x\in[1,2].\end{array}\right.$$
%
\item We have 
\begin{equation}\label{eq:theta-upper-bound}
|\hat\theta_{n,b} - \theta| \le 2\normsup{\hat{g}_b^{sym} - g^{sym}}{[1-\delta,1+\delta]}.
\end{equation}
\item Moreover, 
\begin{equation}\label{eq:markov}
\mathbb{P}\left(\tilde\theta_{n,b}\neq\hat\theta_{n,b}\right)\leq \frac{4}{\delta^2}\E\left[ |\hat\theta_{n,b} - \theta|^2 \right],
\end{equation}
and there exists a constant $C>0$, which only depends on $\delta$, such that
\begin{equation}\label{eq:risque_theta_hat}
\E \left[|\tilde{\theta}_{n,b} - \theta|^2 \right] \leq C\E \left[ \normsup{\hat{g}_b^{sym} - g^{sym}}{[1-\delta,1+\delta]}^2\right]. 
\end{equation}

\end{itemize}
\end{lemma}


The first property of Lemma \ref{lem:g_sym} permits to deal with  $\hat g^{sym}_b$ as with a classical kernel density estimate defined from an \textit{i.i.d} sample. Thus we have $\mathbb E[\hat g_b^{sym}(x)]=L_b\star g^{sym}(x)$. This permits to obtain an upper-bound for the risk of $\hat g_b^{sym}$ as an estimator of $g^{sym}$, and also to define an automatic bandwidth selection rule like for classical kernel density estimates (see paragraph just below).  The second property \eqref{eq:theta-upper-bound} allows us to control the estimation risk of $\hat{\theta}_{n,b}$, while the third one, \eqref{eq:markov}, justifies that the introduction of $\tilde\theta_{n,b}$ is reasonable. 

To obtain a fully data-driven estimate  $\tilde\theta_{n,b}$, it remains to define a bandwidth selection rule for the (classical) kernel estimator $\hat{g}^{sym}_{b}$.  In view of  \eqref{eq:theta-upper-bound}, we introduce a data-driven procedure under sup-norm loss, inspired from Lepski \cite{Lepski2013_sup-norm-loss}. For any $x\in [0,2]$ and any bandwidth $b,b'$ in a collection $\mathcal B'$, we set $\hat g_{b,b'}^{sym}(x) = (L_b\star \hat{g}^{sym}_{b'})(x)$, and $\Gamma_2(b) = \lambda \normsupglobal{L}\log(n)/(nb)$, with $\lambda$ a tuning parameter. As for the other bandwidth selection device, we now define
\[
\Delta(b) = \underset{b' \in \mathcal B'}{\max}\left\{ \underset{x\in [1-\delta, 1+\delta]}{\sup} \big(  \hat g_{b,b'}^{sym}(x) - \hat{g}^{sym}_{b'}(x)\big)^2 - \Gamma_2(b') \right\}_+,
\]
Finally, we choose $\tilde{b} ={\argmin}_{b\in \mathcal B'} \{ \Delta(b) + \Gamma_2(b)\}$, which leads to $\hat{g}^{sym}:= \hat{g}_{\tilde b}^{sym}$ and $\tilde\theta_{n} :=  \tilde\theta_{n,\tilde b}.$ The results of   \cite{Lepski2013_sup-norm-loss} prove that $\mathbb E[\normsup{\hat{g}_b^{sym} - g^{sym}}{[1-\delta,1+\delta]}^2]\leq C(\log n/n)^{2\beta/(2\beta + 1)}$, if $g\in\Sigma(\beta,\mathcal L')$, where $C,\mathcal L'>0$ are some constants, and if the kernel $L$ has an order $\ell =  \lfloor \beta \rfloor$. Combined with Lemma \ref{lem:g_sym}, this ensures that $\tilde\theta_{n} $ satisfies 
\begin{equation}\label{eq:rate-cv-theta}
\E  \left[|\tilde\theta_n - \theta |^2 \right] \leq C\left( \frac{\log n}{n}\right)^{\frac{2\beta}{2\beta + 1}}.
\end{equation}
Numerical simulations in Section \ref{sec:simus} justify that our estimator has a good performance from the practical point of view, in comparison with those proposed in \cite{Ng-Matias-SJS} and \cite{storey2002}.

\subsection{Convergence rate of the component density estimator}\label{sec:cv_rate}

We have now everything we need to compute the convergence rate of our estimator $\hat f(x_0)$ at the point $x_0$, with selected bandwidth, and defined with the preliminary estimates $\hat g$ and $\tilde\theta_{n}$ introduced above (sections \ref{sec:estimation-theta} and \ref{sec:estimation-g} respectively). Starting from the results of Theorem \ref{th:oracle-inequality}, we obtain the following rate of decrease for the pointwise risk of our estimate, over H\"older smoothness classes.

\begin{corollary}\label{cor:rate-cv-fhat}
Assume that {\bf (A1)}, {\bf (A6)} and {\bf (A7)} are satisfied,  for $\beta >0$ and $\LL>0$, and for an index $\ell >0$ such that $\ell \ge \lfloor \beta \rfloor$. We choose e.g. $\alpha_n=\log(n)$ and the bandwidth collection
\[
\mathcal H_n:=\left\{\frac1k, k\in\{1,\hdots,\lfloor\sqrt{n}\rfloor\}\cap[\alpha_n,\widehat\gamma n/\log^3(n)]\right\}.
\]

Then, if $\hat f$ is defined with the preliminary estimates $\hat g$ and $\tilde\theta_{n}$ introduced in sections \ref{sec:estimation-theta} and \ref{sec:estimation-g} respectively, it satisfies
\begin{equation}\label{eq:cv}
 \E\left[\big(  \hat f(x_0) - f(x_0)\big)^2 \right] \le C_{9}^*\left( \frac{\log n}{n}\right)^{\frac{2\beta}{2\beta + 1}},
\end{equation}
where $C_{9}^*$ is a constant depending on $\normsup{g}{\mathcal{V}_n(x_0)}$, $\norme{K}{1}$, $\norme{K}{2}$, $\mathcal{L}$ and $\normsup{f}{\mathcal{V}_n(x_0)}$.
\end{corollary}

The estimator $\hat f$,  with data-driven bandwidth, now achieves the convergence rate $(\log n/n)^{{2\beta}/(2\beta + 1)}$ over the class $\Sigma(\beta, \LL)$ as soon as $\beta \le \ell$. The risk decreases at the optimal minimax rate of convergence (up to a logarithmic term) : the upper bound of Corollary \ref{cor:rate-cv-fhat} matches with the lower-bound for the minimax risk established by Ibragimov and Hasminskii \cite{ibragimov_estimate_1980}. 
Our procedure automatically adapts to the unknown smoothness of the function to estimate : the bandwidth $\hat h(x_0)$ is computed in a fully data-driven way, without using the knowledge of the regularity index $\beta$, contrary to the estimator $\hat{f}_n^{rwk}$ of Nguyen and Matias \cite{nguyen_matias_2014} (corollary 3.4). 

\begin{remark}
In the present work, we focus on  Model \eqref{model}. However, the estimation procedure we develop  can easily  be extended to the model 
\begin{equation}\label{general-model}
g(x) = \theta\phi(x) + (1-\theta)f(x), \quad x\in \R,
\end{equation}
where the function $\phi$ is a known density, but not necessarily equal to the uniform one.  In this case, a family of kernel estimates can be defined  like in \eqref{eq:estimator-f} replacing the weights $w(\tilde\theta_n, \hat g(\cdot))$ by 
\[
w(\tilde \theta_n, \hat{g}(\cdot), \phi(x_0)) = \frac 1{1-\tilde \theta_n}\left(1 - \frac{\tilde \theta_n\phi(x_0)}{\hat{g}(\cdot)}\right).
\]
If the density function $\phi$ is uniformly bounded on $\R$, it is then possible to obtain analogous results (bias-variance trade-off for the pointwise risk,  adaptive bandwidth selection rule leading to oracle-type inequality and optimal convergence rate) as we established for model \eqref{model}. 
\end{remark}

\section{Numerical study}\label{sec:simus}

\subsection{Simulated data}\label{sec:simu_model}

\begin{figure}[!ht]\begin{center}
\begin{tabular}{ccc}
\includegraphics[scale=0.28]{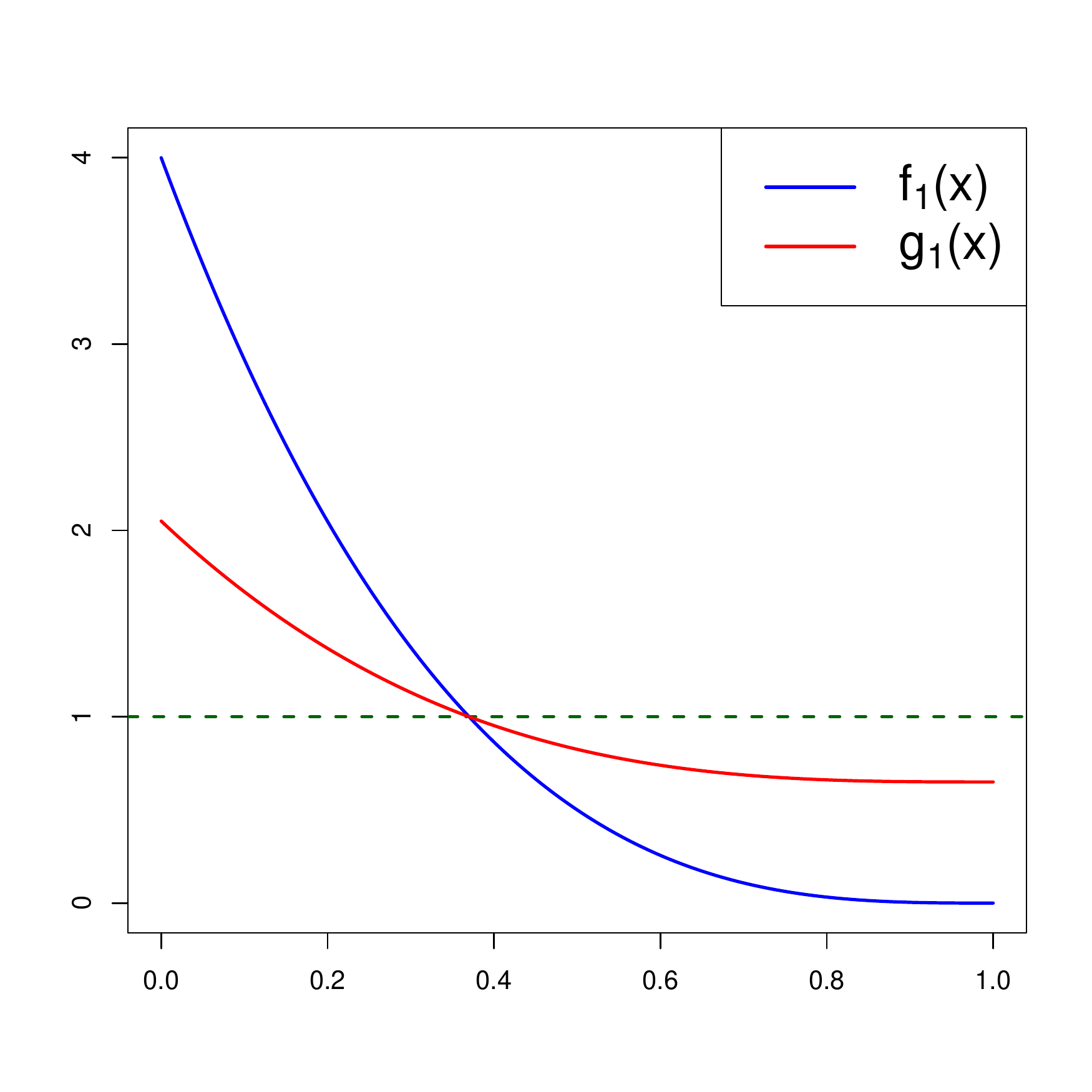} &
\includegraphics[scale=0.28]{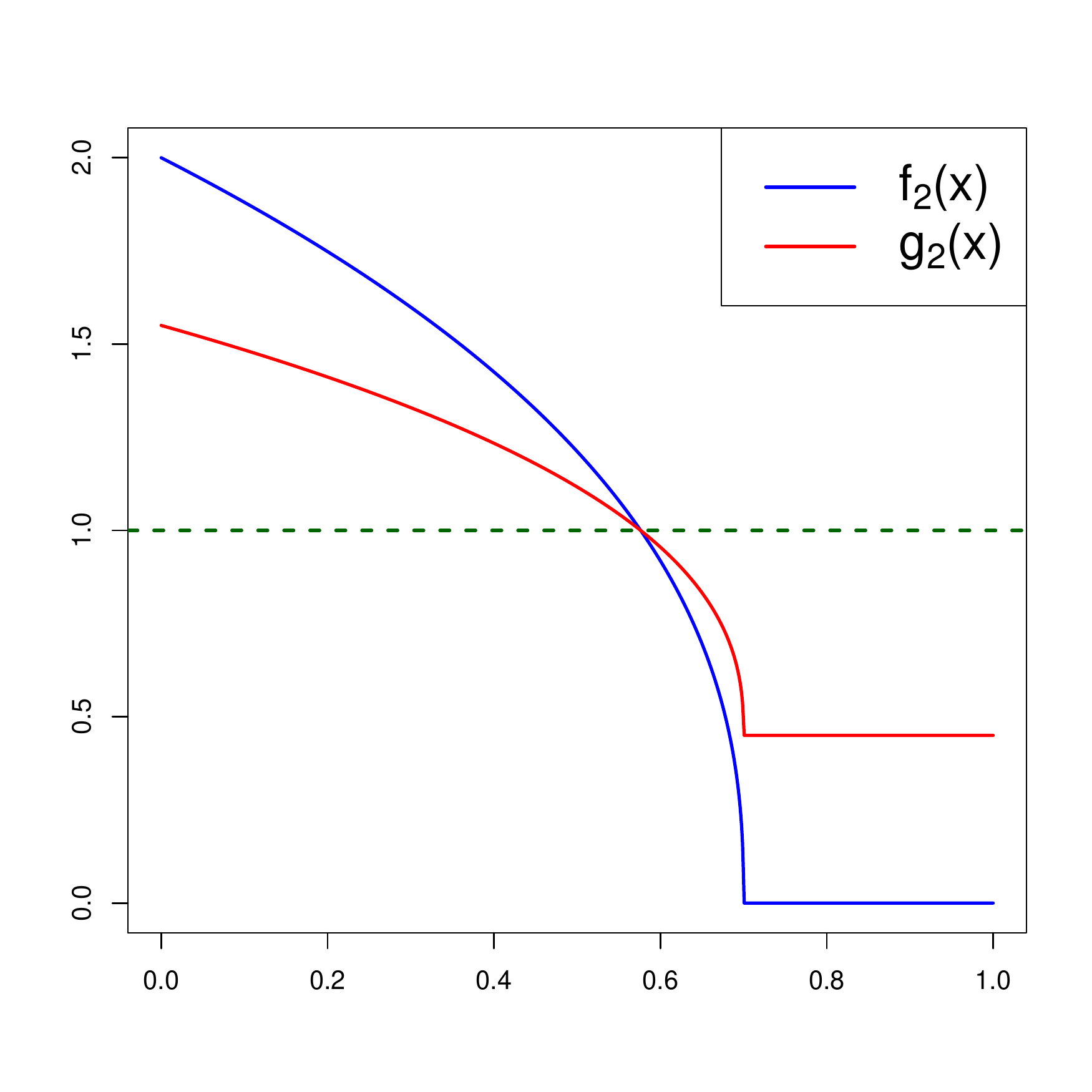} &
\includegraphics[scale=0.28]{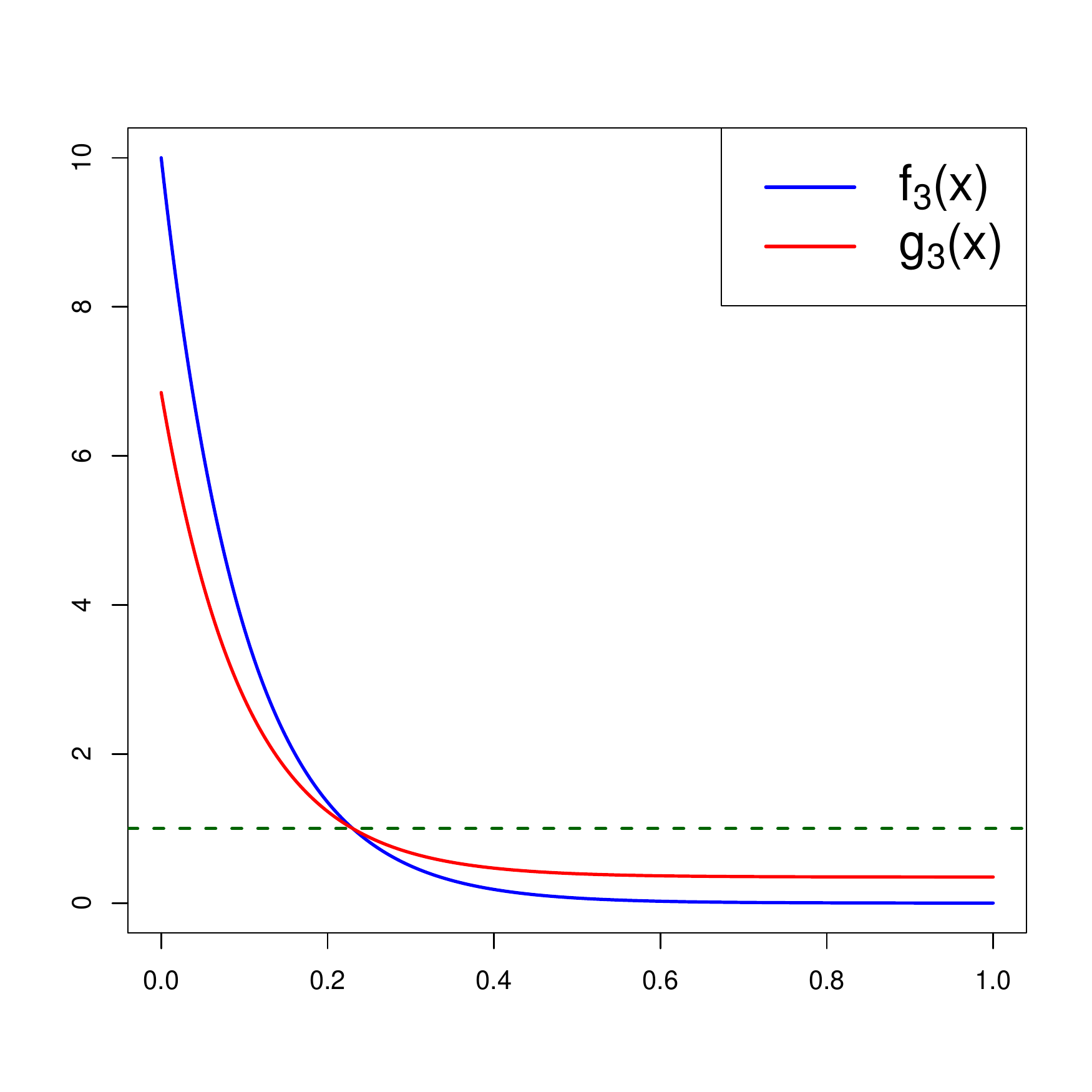} 
\end{tabular}
\caption{Representation of $f_j$ and the corresponding $g_j$ in model \eqref{model} for $(\theta_1 = 0.65,f_1)$ (left), $(\theta_2 = 0.45, f_2)$ (middle) and $(\theta_3 = 0.35, f_3)$ (right).}
\label{fig:design-densities}
\end{center}\end{figure}

We briefly  illustrate the performance of the estimation method over simulated data, according the following framework.  We simulate observations with density $g$ defined by model \eqref{model} for sample size $n\in \{500, 1000, 2000\}$. Three different cases of $(\theta,f)$ are considered: \begin{itemize}
\item $f_1(x) = 4(1-x)^3\id_{[0,1]}(x)$, $\theta_1 = 0.65$. 

\item $f_2(x) = \dfrac{s}{1-\delta}\left(1 - \dfrac{x}{1-\delta}\right)^{s-1}\id_{[0,1-\delta]}(x)$ with $(\delta, s) = (0.3, 1.4)$, $\theta_2 = 0.45$. 

\item $f_3(x) = \lambda e^{-\lambda x}\left( 1 - e^{-\lambda b} \right)^{-1} \id_{[0,b]}(x)$ the density of truncated exponential distribution on $[0,b]$ with $(\lambda, b) = (10, 0.9)$, $\theta_3 = 0.35$.  
\end{itemize}
The density $f_1$ is borrowed from \cite{nguyen_matias_2014} while the shape of $f_2$ is used both by \cite{celisse2010} and \cite{Ng-Matias-SJS}. Figure \ref{fig:design-densities} represents those three cases with respect to each design density and associated proportion $\theta$.

\subsection{Implementation of the method}\label{sec:tuning}

To compute our estimates, we choose $K(x) = L(x) = (1 - |x|)\id_{\{|x| \le 1\}}$ the triangular kernel. In the variance term \eqref{eq:kappa-variance-term} of the GL method used to select the bandwidth of the kernel estimator of $f$, we replace  $\|g\|_{\infty,\mathcal{V}_n(x_0)}$ by the $95^{\text{th}}$ percentile of $\big\{\max_{t\in \mathcal{V}_n(x_0)} \hat g_h(t), h\in \Hn\big\}$. Similarly, in the variance term $\Gamma_1$ used to select the bandwidth of the kernel estimate of $g$, we use the $95^{\text{th}}$ percentile of $\big\{\max_{t\in [0,1]} \hat g_h(t), h\in \Hn\big\}$ instead of $\|g\|_{\infty}$. The collection of bandwidths $\mathcal H_n, \mathcal B, \mathcal B'$ are equal to $\big\{1/k, k = 1,\ldots, \lfloor \sqrt{n} \rfloor \big\}$.

\begin{figure}[htbp]
\centering
\begin{tabular}{ccc}
\includegraphics[scale=0.28]{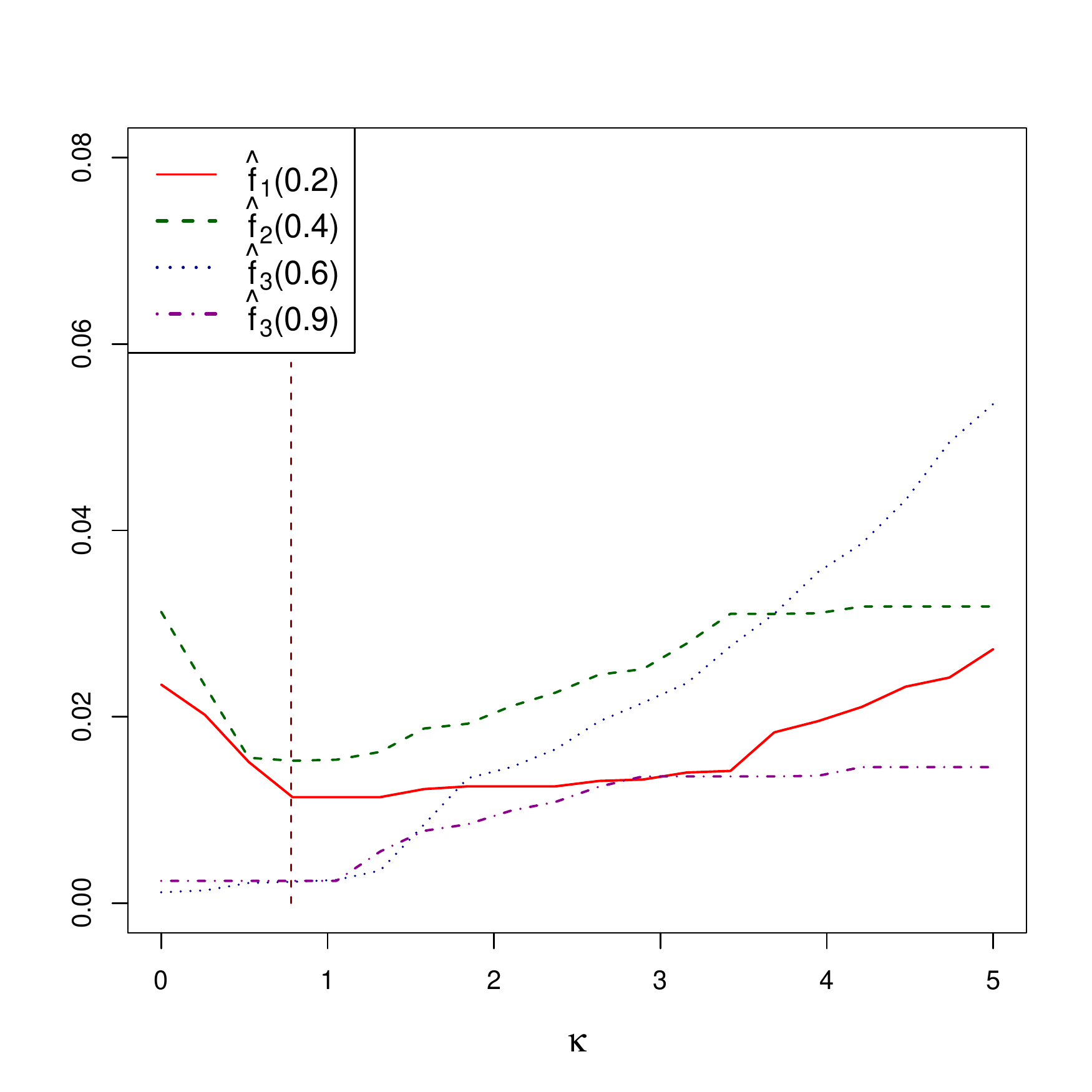} &  \includegraphics[scale=0.28]{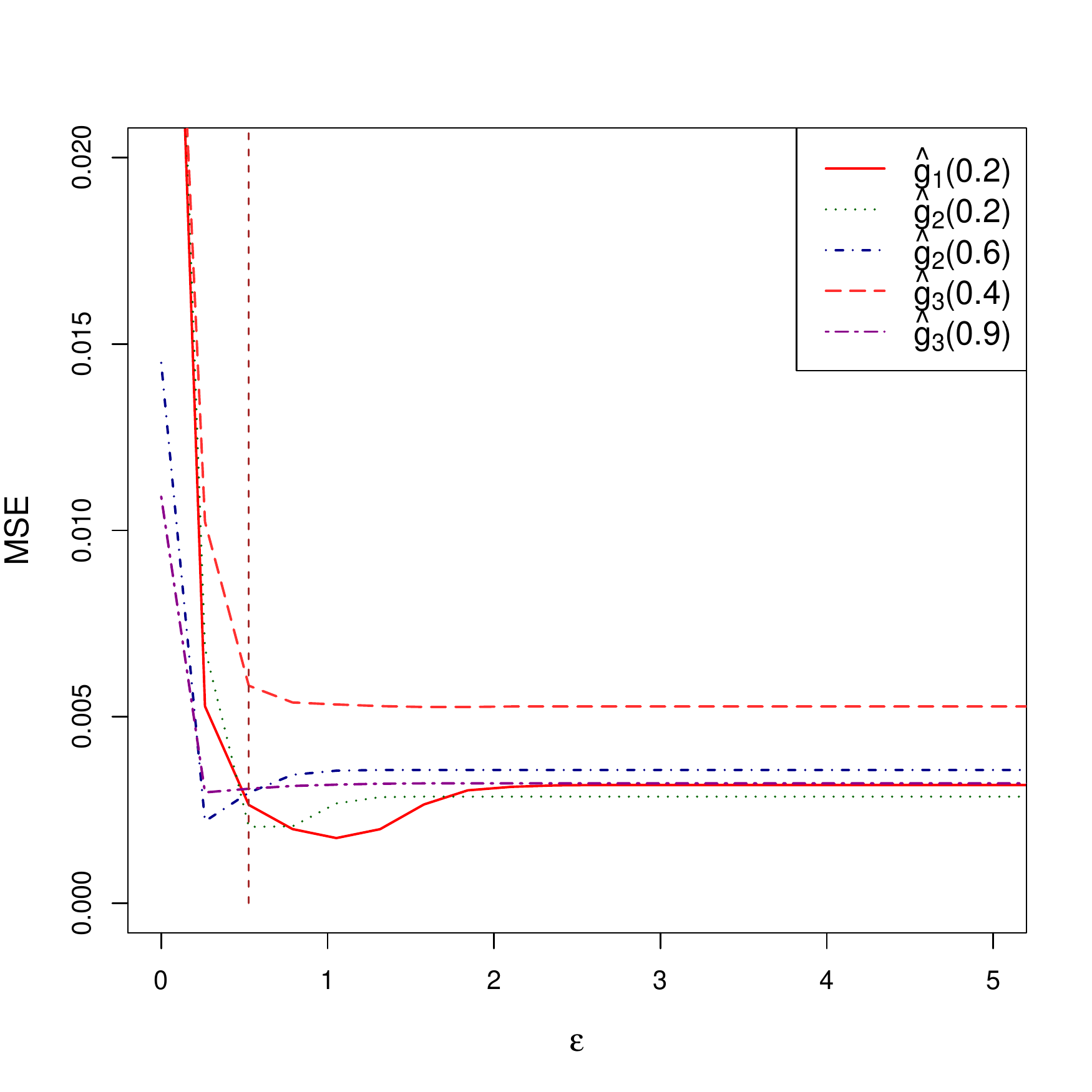}   & \includegraphics[scale=0.28]{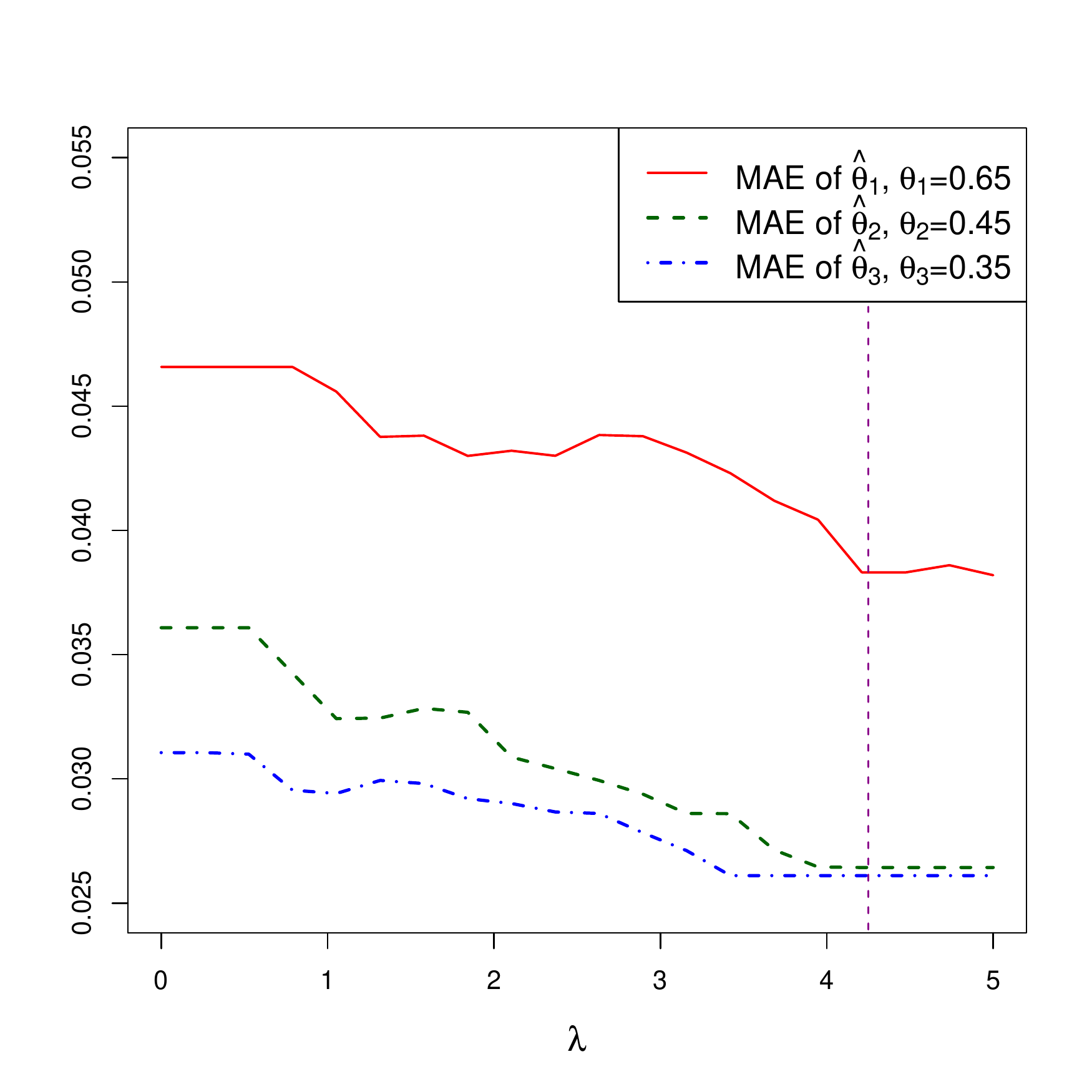} \\
(a) & (b) & (c) \\
\end{tabular}
\caption{values of the mean-squared error for (a)  $\hat{f}(x_0)$ with respect to $\kappa$, (b) $\hat{g}(x_0)$ with respect to $\varepsilon$.  (c) : Values of the mean-absolute error for $\hat\theta_n$  with respect to $\lambda$. The sample size is $n=2000$ for all computations. The vertical line corresponds to the chosen value of $\kappa$ (figure (a)), $\varepsilon$ (figure (b)) and $\lambda$ (figure (c)).}
\label{fig:calibration-constants}
\end{figure}

We shall  settle the values of the constants $\kappa$, $\varepsilon$ and $\lambda$ involved in the penalty terms $V(x_0,h), \Gamma_1(h)$ and $\Gamma_2(b)$ respectively, to compute the selected bandwidths. Since the calibrations of these tuning parameters are carried out in the same fashion, we only describe the calibration for $\kappa$.  Denote by $\hat{f}_{\kappa}$ the estimator of $f$ depending on the constant $\kappa$ to be calibrated. We approximate the mean-squared error for the estimator $\hat{f}_{\kappa}$, defined by $\MSE(\hat f_{\kappa}(x_0)) = \E[( \hat f_{\kappa}(x_0) - f(x_0) )^2 ]$, over 100 Monte-Carlo runs, for different possible values $\{\kappa_1, \ldots, \kappa_K \}$ of  $\kappa$, for the three densities $f_1$, $f_2$, $f_3$ calculated at several test points $x_0$. We choose a value for $\kappa$ that leads to small risks in all investigated cases.   Figure \ref{fig:calibration-constants}\red{(a)} shows that $\kappa = 0.78$ is an acceptable choice even if other values can be also convenient. Similarly, we  set $\varepsilon = 0.52$ and $\lambda = 4.25$ (see Figure \ref{fig:calibration-constants}\red{(b)}  and \ref{fig:calibration-constants}\red{(c)} for the calibrations of $\varepsilon$ and $\lambda$. 

\subsection{Simulation results}

\subsubsection{Estimation of the mixing proportion $\theta$} 

\begin{figure}[!ht]
\begin{center}
\begin{tabular}{ccc}
\includegraphics[scale=0.28]{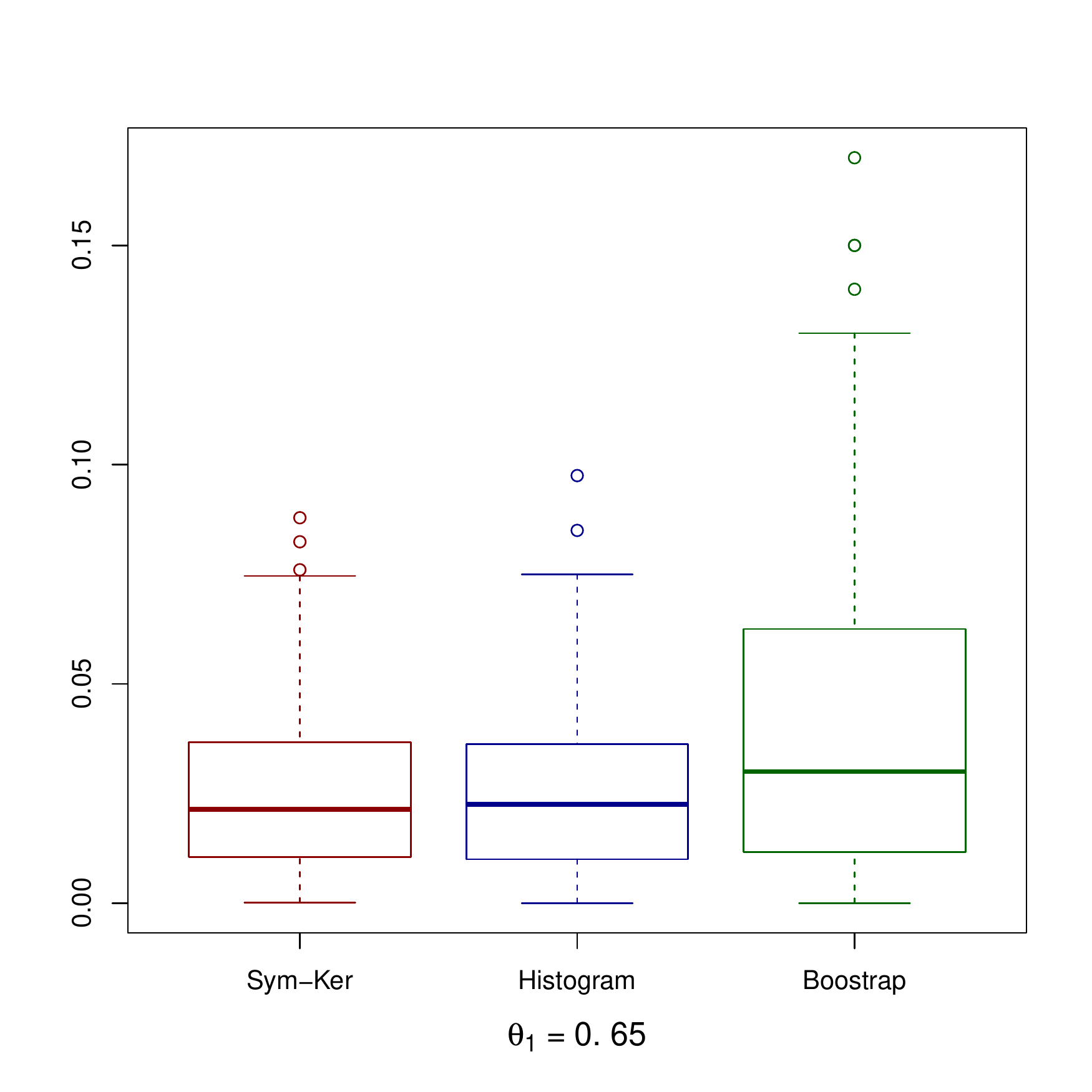} & \includegraphics[scale=0.28]{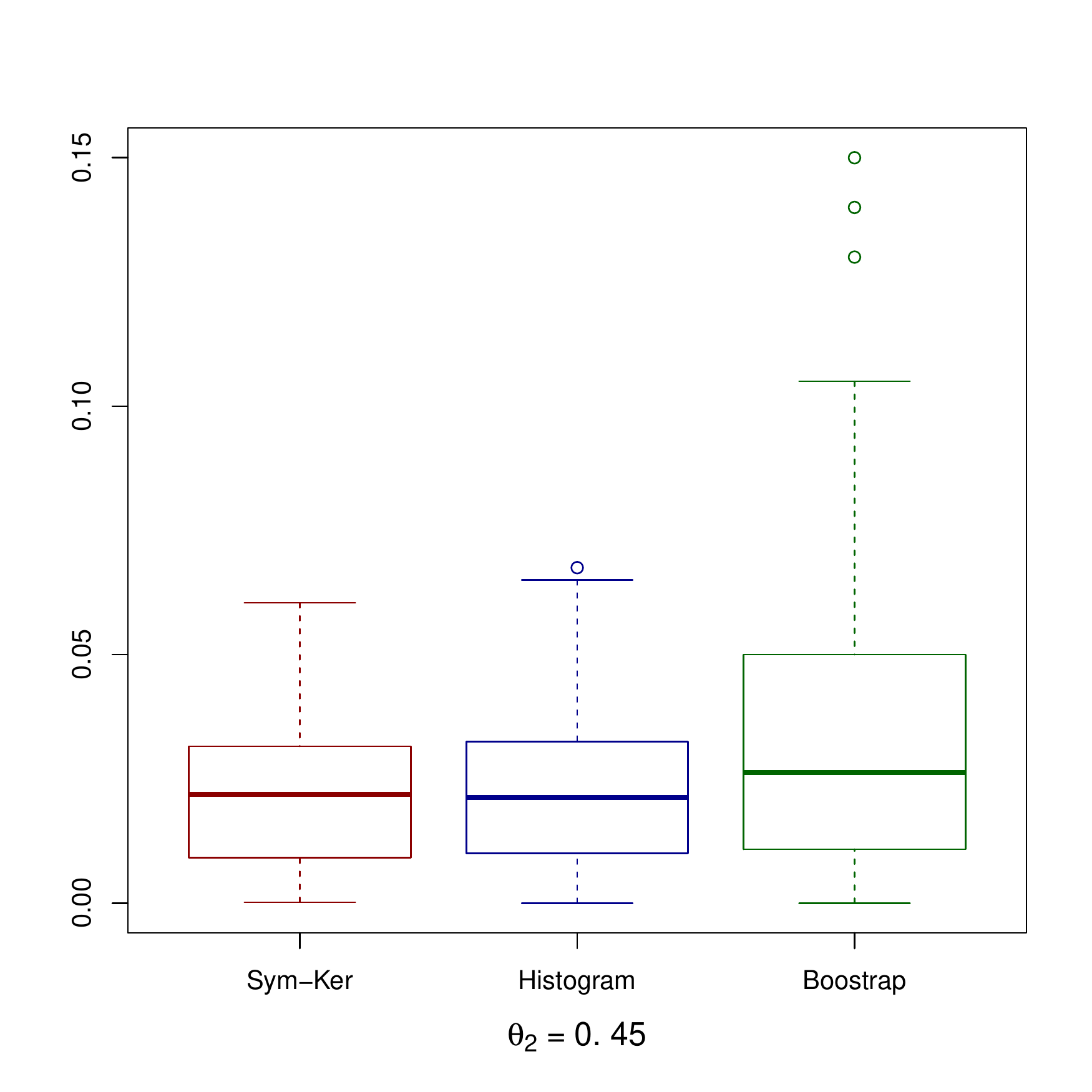} &\includegraphics[scale=0.28]{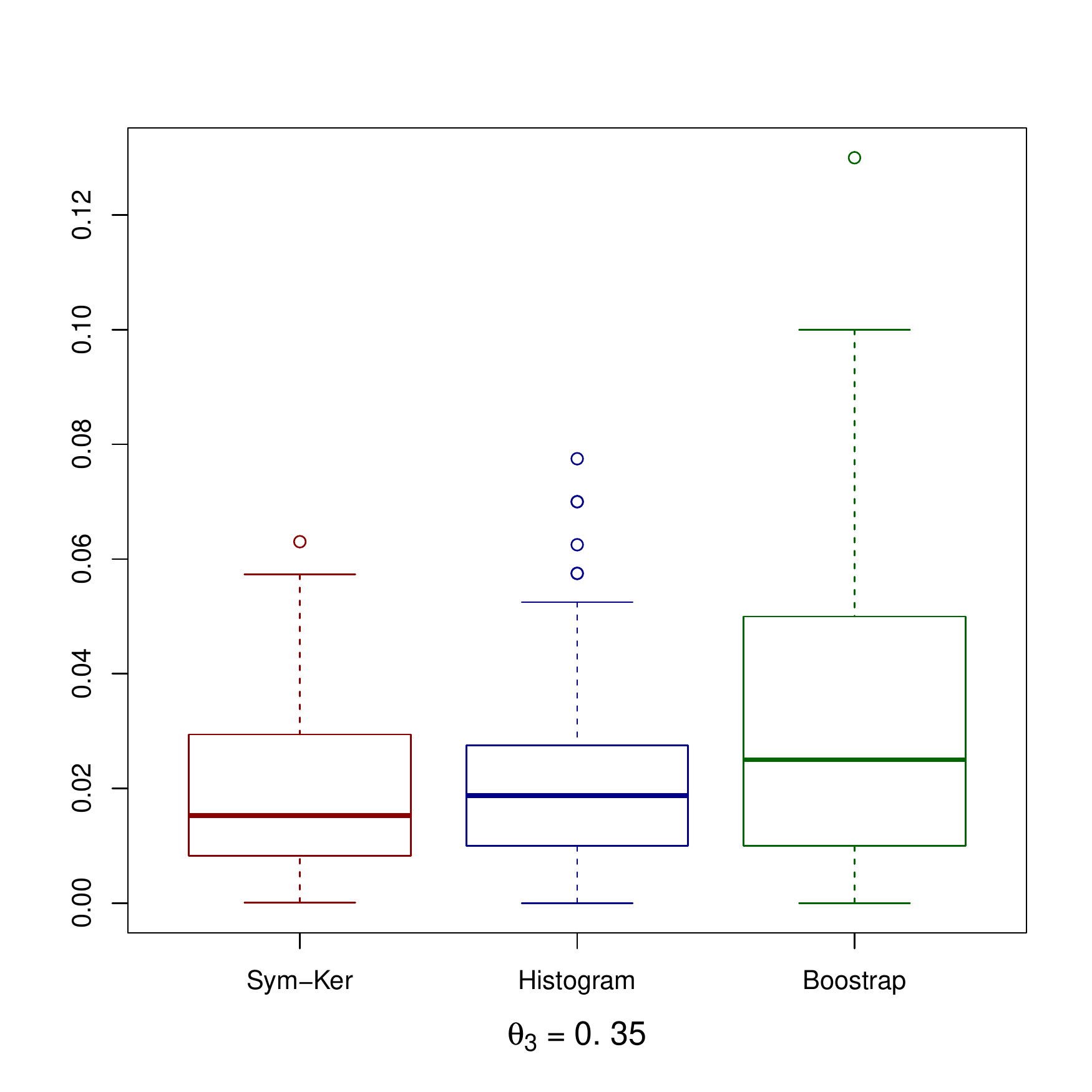}  \\
\end{tabular}
\caption{errors for the estimation of $\theta$ in the three simulated settings (with sample size $n=2000$).}
\label{fig:err-estimate-theta}
\end{center}
\end{figure}

We compare our estimator $\hat\theta_n$ with the histogram-based estimator $\hat\theta_n^{\text{Ng-M}}$ proposed in \cite{Ng-Matias-SJS} and the estimator $\hat\theta_n^{S}$ introduced in \cite{storey2002}. Boxplots in Figure \ref{fig:err-estimate-theta} represent the absolute errors of $\hat\theta_n$, $\hat\theta_n^{\text{Ng-M}}$ and $\hat\theta_n^{S}$, labeled respectively by "Sym-Ker", "Histogram" and "Bootstrap". The estimators $\hat\theta_n$  and  $\hat\theta_n^{\text{Ng-M}}$ have comparable  performances, and outperform $\hat\theta_n^{S}$.  

\subsubsection{Estimation of the target density $f$} 

\begin{figure}[!ht]
\centering
\includegraphics[scale = 0.7]{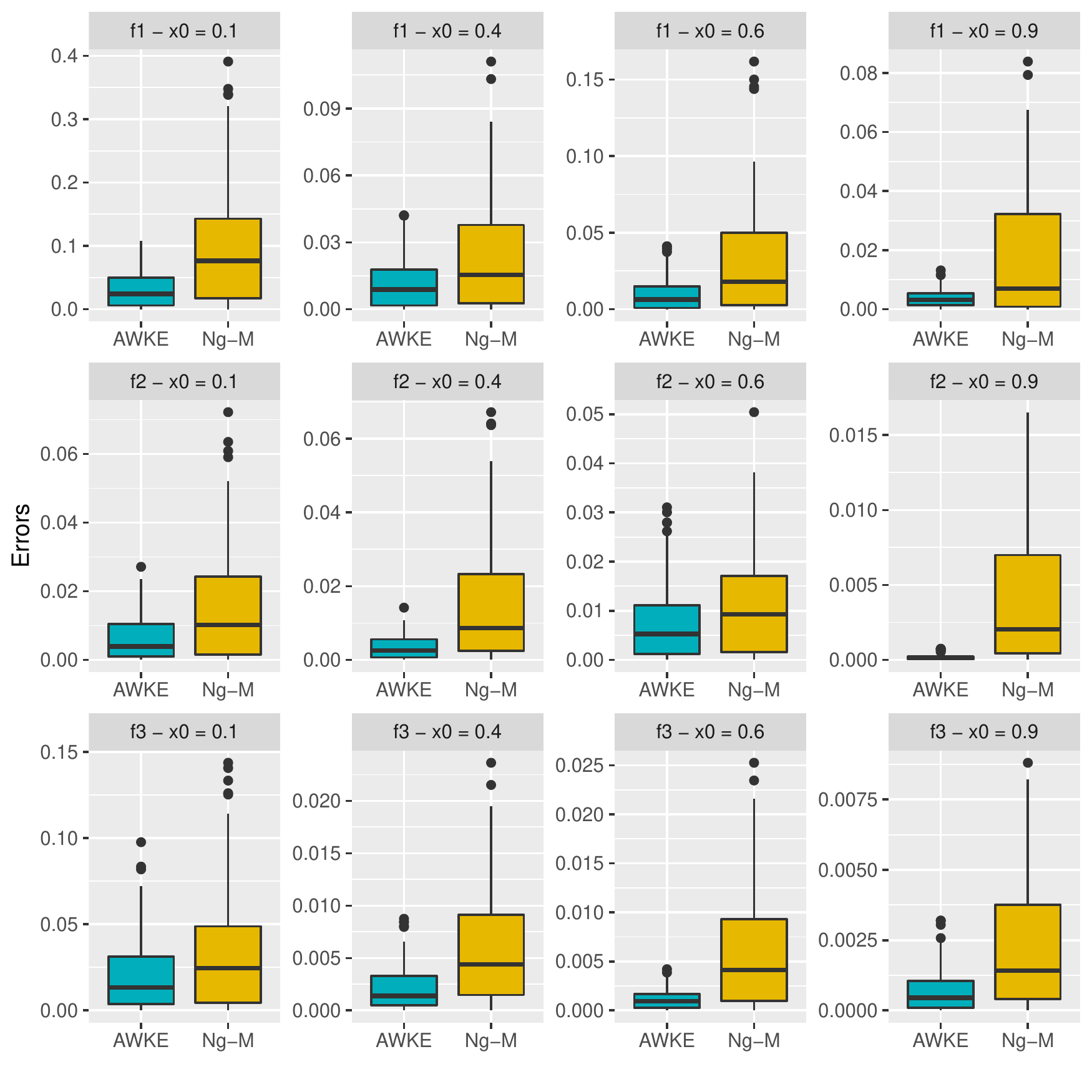}
\caption{errors for the estimation of $f_1$, $f_2$ and $f_3$ for $x_0 \in \{0.1, 0.4, 0.6, 0.9\}$ and sample size $n=2000$. \label{fig:err-fhat-n2000}}
\end{figure}

We present in Tables \ref{tab:MSE-fhat1}, \ref{tab:MSE-fhat2} and \ref{tab:MSE-fhat3} the mean-squared error (MSE) for the estimation of $f$ according to the three different models and the different sample sizes  introduced in Section \ref{sec:simu_model}. The MSEs' are approximated over $100$ Monte-Carlo replications. We shall choose the estimation points (to compute the pointwise risk): we propose $x_0 \in \{0.1, 0.4, 0.6, 0.9\}$.  The choices of $x_0 = 0.4$ and $x_0 = 0.6$ are standard. The choices of $x_0 = 0.1$ and $x_0 = 0.9$ allows to test the performance of $\hat{f}$ close to the boundaries of the domain of definition of $f$ and $g$.  We compare our estimator $\hat f$ with the randomly weighted estimator proposed in Nguyen and Matias \cite{nguyen_matias_2014}. In the sequel, the label ''AWKE'' (Adaptive Weighted Kernel Estimator) refers to our estimator $\hat f$, whose bandwidth is selected by the Goldenshluger-Lepski method and ''Ng-M'' refers to the one proposed by \cite{nguyen_matias_2014}.  Resulting boxplots are displayed in  Figure \ref{fig:err-fhat-n2000} for $n=2000$. 

\renewcommand{\arraystretch}{1.15}
\begin{table}[htp]
\begin{center}
\begin{tabular}{llcccc}
\hline
Sample size 	& Estimator 	& $x_0 = 0.1$ 	& $x_0 = 0.4$ 	& $x_0 = 0.6$ 	& $x_0 = 0.9$ \\
\hline
$n = 500$ 	&AWKE		& 0.1683		& 0.0119		& 0.0256		& 0.0059 \\
			&Ng-M 		& 0.2869 		& 0.0450		& 0.1046		& 0.0433 \\
\hline
$n = 1000$	&AWKE		& 0.0632		& 0.0087		& 0.0118		& 0.0063 \\
			&Ng-M		& 0.1643		& 0.0469		& 0.0651		& 0.0279 \\
\hline 
$n = 2000$	&AWKE		& 0.0314		& 0.0118		& 0.0098		& 0.0038 \\
			&Ng-M		& 0.0982		& 0.0246		& 0.0326		& 0.0164 \\
\hline 
\end{tabular}
\caption{mean-squared error of  the reconstruction of $f_1$, for our estimator $\hat f$ (AWKE), and for the estimator of Nguyen and Matias \cite{nguyen_matias_2014} (Ng-M).  \label{tab:MSE-fhat1}}
\end{center}
\end{table}%
    
\begin{table}[!ht]
\begin{center}
\begin{tabular}{llcccc}
\hline
Sample size 	& Estimator 	& $x_0 = 0.1$ 	& $x_0 = 0.4$ 	& $x_0 = 0.6$ 	& $x_0 = 0.9$ \\
\hline
$n = 500$ 	&AWKE 		& 0.0430 	& 0.0126 		& 0.0311		& 0.0002 \\
			&Ng-M 		& 0.0560 		& 0.0540		& 0.0306		& 0.0138 \\
\hline
$n = 1000$	&AWKE		& 0.0183		& 0.0061		& 0.0240		& 0.0005 \\
			&Ng-M		& 0.0277		& 0.0209		& 0.0123		& 0.0069 \\
\hline 
$n = 2000$	&AWKE		& 0.0061		& 0.0034 	& 0.0076	& 0.0002 \\
			&Ng-M		& 0.0164		& 0.0159		& 0.0113		& 0.0038 \\
\hline 
\end{tabular}
\caption{mean-squared error of  the reconstruction of $f_2$, for our estimator $\hat f$ (AWKE), and for the estimator of Nguyen and Matias \cite{nguyen_matias_2014} (Ng-M). \label{tab:MSE-fhat2}}
\end{center}
\end{table}%
         
\begin{table}[!ht]
\begin{center}
\begin{tabular}{llcccc}
\hline
Sample size 	& Estimator 	& $x_0 = 0.1$ 	& $x_0 = 0.4$ 	& $x_0 = 0.6$ 	& $x_0 = 0.9$ \\
\hline
$n = 500$ 	&AWKE		& 0.0737 		& 0.0090 	& 0.0039		& 0.0016 \\
			&Ng-M 		& 0.1308 		& 0.0247		& 0.0207		& 0.0096 \\
\hline
$n = 1000$	&AWKE		& 0.0296		& 0.0051		& 0.0026		& 0.0009 \\
			&Ng-M		& 0.0566		& 0.0106		& 0.0096		& 0.0060 \\
\hline 
$n = 2000$	&AWKE		& 0.0224		& 0.0022	 	& 0.0012		& 0.0007 \\
			&Ng-M		& 0.0342		& 0.0059		& 0.0062		& 0.0021 \\
\hline 
\end{tabular}
\caption{mean-squared error of  the reconstruction of $f_3$, for our estimator $\hat f$ (AWKE), and for the estimator of Nguyen and Matias \cite{nguyen_matias_2014} (Ng-M). \label{tab:MSE-fhat3}}
\end{center}
\end{table}%

Tables \ref{tab:MSE-fhat1}, \ref{tab:MSE-fhat2}, \ref{tab:MSE-fhat3} and boxplots show that our estimator outperforms the one of \cite{nguyen_matias_2014}. Notice that the errors are relatively large at the point $x_0 = 0.1$, for both estimators, which was expected (boundary effect). 

\newpage
\section{Proofs}\label{sec:proofs}
In the sequel, the notations $\tilde\Prob$, $\tilde\E$ and $\tilde{\mathbb{V}}ar$ respectively denote the probability, the expectation and the variance associated with $X_1,\ldots, X_n$, conditionally on the additional random sample $X_{n+1}, \ldots, X_{2n}$. 

\subsection{Proof of Proposition \ref{prop:upper-bound-fhat}}\label{proof:prop-upper-bound-fhat}

Let $\rho > 1$, introduce the event
\[
\Omega_\rho = \left\{ \rho^{-1}\gamma \le \hat\gamma \le \rho\gamma \right\}.
\]
such that
\begin{equation}\label{eq:risk-decompose}
\hat f_h(x_0) - f(x_0)=  \big(\hat f_h(x_0) - f(x_0)\big) \id_{\Omega_\rho} + \big(\hat f_h(x_0) - f(x_0)\big) \id_{\Omega^c_\rho}.
\end{equation}

We first evaluate the term $\big( \hat f_h(x_0) - f(x_0)\big) \id_{\Omega_\rho}$. Suppose now that we are on $\Omega_\rho$, then for any $x_0\in [0,1]$, we have
\begin{equation}\label{eq:risk-on-Omega}
\big( \hat f_h(x_0) - f(x_0) \big)^2  \le 3\Big( \big( \hat f_h(x_0) - K_h\star \check f(x_0) \big)^2  + \big(  K_h\star \check f(x_0) - \check f(x_0) \big)^2  + 
\big( \check f(x_0) - f(x_0) \big)^2 \Big),
\end{equation}
where we define 
\[
\check f(x) = w(\tilde\theta_n, \hat g(x))g(x) = \frac{1}{1-\tilde\theta_n}\left(1 - \frac{\tilde\theta_n}{\hat g(x)} \right)g(x).
\]
Note that by definition of $\check f$, we have $K_h\star \check f(x_0) = \tilde\E\big[\hat f_h(x_0)\big]$. Hence,
\[
\big( \hat f_h(x_0) - K_h\star \check f(x_0) \big)^2  = \big( \hat f_h(x_0) - \tilde\E\big[\hat f_h(x_0)\big]\big)^2.
\]
It follows that
\begin{align}
\tilde \E \left[\big( \hat f_h(x_0) - \tilde\E\big[\hat f_h(x_0)\big]\big)^2 \right]&= \tilde{\mathbb{V}}ar \left(\hat f_h(x_0) \right) = \tilde{\mathbb{V}}ar \left( \frac 1n\sumi w(\tilde\theta_n, \hat g(X_i))K_h(x_0 - X_i) \right) \notag \\
&=\frac 1n \tilde{\mathbb{V}}ar  \left( w(\tilde\theta_n, \hat g(X_1))K_h(x_0 - X_1)\right) \notag \\
&\le \frac 1n \tilde\E\left[ \left(w(\tilde \theta_n, \hat g(X_1))K_h(x_0 - X_1) \right)^2 \right]. \notag 
\end{align}
On the other hand, for all $i\in \{1,\ldots, n\}$, thanks to {\bf (A4)} and {\bf (A2)}, and since $\hat\gamma \ge \rho^{-1}\gamma$ on $\Omega_\rho$,
\begin{align}
\left| w(\tilde \theta_n, \hat g(X_i))K_h(x_0 - X_i) \right| &= \left| \frac{1}{1 - \tilde\theta_n}\left(1 - \frac{\tilde\theta_n}{\hat g(X_i)} \right)K_h(x_0 - X_i) \right| \le 
\frac{2}{\delta}\left(1 + \frac{\tilde\theta_n}{|\hat g(X_i)|} \right) | K_h(x_0 - X_i) | \notag \\
&\le \frac{2}{\delta}\left(1 + \frac{1}{\hat \gamma} \right) |K_h(x_0 - X_i)| \label{eq:upperbound-w-theta-g-bis}  \\
&\le \frac{2}{\delta}\left(1 + \frac{\rho}{ \gamma} \right) |K_h(x_0 - X_i)| = \frac{2(\rho + \gamma)}{ \delta \gamma} |K_h(x_0 - X_i)|. \label{eq:upperbound-w-theta-g}
\end{align}
For \eqref{eq:upperbound-w-theta-g-bis}, as we use compactly supported kernel to construct the estimator $\hat{f}_h$, condition  $\alpha_n \le h^{-1}$ in {\bf (A5)} ensures that $\left| \left(\hat g(X_i)\right)^{-1}K_h(x_0 - X_i) \right|$ is upper bounded by $\hat\gamma^{-1}|K_h(x_0-X_i)|$ even though we have no observation in the neighbourhood of $x_0$. 

Thus we obtain, on $\Omega_{\rho}$, 

\begin{equation}\label{eq:upper-bound-variance-term}
\tilde \E \left[\big( \hat f_h(x_0) - \tilde\E\big[\hat f_h(x_0)\big]\big)^2 \right] \le  \frac{4(\rho + \gamma)^2}{\delta^{2}\gamma^{2}n}  \tilde\E\left[K_h^2(x_0 - X_1) \right]\le
 \frac{4 (\rho + \gamma)^2 \|K\|_2^2\normsup{g}{\mathcal{V}_n(x_0)}}{\delta^{2}\gamma^{2}nh}. 
\end{equation}

For the last two terms of \eqref{eq:risk-on-Omega}, we apply the following proposition, which proof can be found in Section \ref{sec:proof-prop-control-bias}.

\begin{proposition}\label{prop:control-bias-fcheck}
Assume {\bf (A1)} and {\bf (A3)}. On the set $\Omega_\rho$, we have the following results for any $x_0\in [0,1]$
\begin{align}
\big(\check f(x_0) - f(x_0) \big)^2 &\le C_1 \delta^{-2}\gamma^{-2} \normsup{\hat g - g}{\mathcal{V}_n(x_0)}^2 + C_2 \delta^{-6}\big|\tilde{\theta}_n - \theta\big|^2, \label{eq:upperbound-fcheck-f}\\
\big(  K_h\star \check f(x_0) - \check f(x_0) \big)^2 &\le 6 \normsup{ K_h\star f - f }{\mathcal{V}_n(x_0)}^2  + C_3 \delta^{-2}\gamma^{-2}  \normsup{\hat g - g}{\mathcal{V}_n(x_0)}^2 + C_4\delta^{-6}\big|\tilde{\theta}_n - \theta\big|^2, \label{eq:control-bias-fcheck}
\end{align}
where $C_1$ and $C_2$ respectively depend on $\rho$ and $\norme{g}{\infty, \mathcal{V}_n(x_0)}$,  $C_3$ depends on $\rho$ and $\norme{K}{1}$ and $C_4$ depends on $\norme{g}{\infty, \mathcal{V}_n(x_0)}$ and $\norme{K}{1}$. 
\end{proposition}

Combining \eqref{eq:upper-bound-variance-term}, \eqref{eq:upperbound-fcheck-f} and \eqref{eq:control-bias-fcheck}, we obtain 
\begin{multline*}
\E\left[ \big( \hat f_h(x_0) - f(x_0) \big)^2 \id_{\Omega_\rho} \right] \le  18 \normsup{ K_h\star f - f }{\mathcal{V}_n(x_0)}^2 + 3(C_1+C_3) \delta^{-2}\gamma^{-2} \E\big[\normsup{\hat g - g}{\mathcal{V}_n(x_0)}^2\big] \\ + 3(C_2+C_4) \delta^{-6}\E\big[ \big|\tilde{\theta}_n - \theta\big|^2\big] + \frac{12(\rho+\gamma)^2\norme{K}{2}^2\normsup{g}{\mathcal{V}_n(x_0)}}{\delta^2\gamma^2 nh}.
\end{multline*}

It remains to study the risk bound on $\Omega^c_\rho$. To do so, we successively apply the following lemmas whose proofs are postponed to Section \ref{sec:intermediate_proof}.

\begin{lemma}\label{lem:control-Omega-rho}
Suppose that Assumption {\bf (A3)} is satisfied. Then we have for $\rho>1$
\[
\Prob\left( \Omega_\rho^c \right) \le C_{g,\rho} \exp\left\{ -(\log n)^{3/2} \right\},
\]
with $C_{g,\rho}$ a positive constant depending on $g$ and $\rho$. 
\end{lemma}

\begin{lemma}\label{lem:control-risk-on-OmegaC}
Assume {\bf (A1)} to {\bf (A5)}. For any $h\in \Hn$, we have
\[
\E\left[ \big( \hat f_h(x_0) - f(x_0)\big)^2 \id_{\Omega^c_\rho}\right] \le \frac{C_4^*}{ n^2},
\]
with $C_4^*$ a positive constant depending on $\normsup{f}{\mathcal{V}_n(x_0)}$, $\normsupglobal{K}$, $g$, $\delta$ and $\rho$. 
\end{lemma}

This concludes the proof of Proposition \ref{prop:upper-bound-fhat}.

\subsection{Proof of Theorem \ref{th:oracle-inequality}}\label{sec:proof_thm}
Suppose that we are on $\Omega_\rho$. Let $\hat f$ be the adaptive estimator defined in \eqref{eq:adaptive-estimator-f}, we have for any $x_0\in [0,1]$,
\[
\big( \hat f(x_0)  - f(x_0) \big)^2 \le 2\left(\big( \hat f(x_0) - \check f(x_0) \big)^2 + \big( \check f(x_0) - f(x_0) \big)^2\right)
\]

The second term is controlled by \eqref{eq:upperbound-fcheck-f} of Proposition \ref{prop:control-bias-fcheck}. Hence it remains to handle with the first term. For any $h\in \Hn$, we have
\begin{align}
\big(\hat f(x_0) - \check f(x_0) \big)^2 &\le 3\left( \big(\hat f_{\hat h(x_0)}(x_0) - \hat f_{\hat h(x_0),h}(x_0) \big)^2 + \big(\hat f_{\hat h(x_0),h}(x_0) - \hat f_h(x_0) \big)^2 + \big(\hat f_h(x_0)  - \check f(x_0) \big)^2\right) \notag \\
&= 3\left(\big(\hat f_{\hat h(x_0)}(x_0) - \hat f_{\hat h(x_0),h}(x_0) \big)^2 - V(x_0, \hat h(x_0)) + \big(\hat f_{\hat h(x_0),h}(x_0) - \hat f_h(x_0) \big)^2 - V(x_0,h) \right. \notag \\
&\left.\hspace{2.3in} + V(x_0,\hat h(x_0)) + V(x_0,h) + \big(\hat f_h(x_0)  - \check f(x_0) \big)^2 \right) \notag \\
&\le 3\left(A(x_0, \hat h(x_0)) + A(x_0, h) + V(x_0,\hat h(x_0)) + V(x_0,h) +  \big(\hat f_h(x_0)  - \check f(x_0) \big)^2 \right) \label{eq:interm_preuve} \\
&\le 6A(x_0, h) + 6V(x_0,h) + 3\big(\hat f_h(x_0)  - K_h \star \check f(x_0) \big)^2 + 3\big(K_h\star \check f(x_0)  - \check f(x_0) \big)^2 . \label{eq:bound-fhat}
\end{align}
To obtain \eqref{eq:interm_preuve}, we use the definition of $A(x_0,h)$, $h\in\mathcal H_n$, see \eqref{eq:def_A}, and also that $\hat{f}_{h,h'}=\hat{f}_{h',h}$, for any $h,h'\in\mathcal H_n$. Then, \eqref{eq:bound-fhat} is a consequence of the definition of $\hat h(x_0)$, see \eqref{eq:adaptive-estimator-f}. Next, we have 
\begin{align*}
A(x_0,h) &= \umax{h'\in \Hn}\,\left\{ \big( \hat f_{h,h'}(x_0) - \hat f_{h'}(x_0) \big)^2 - V(x_0,h') \right\}_+ \\
&\le 3\umax{h'\in \Hn}\left\{ \big( \hat f_{h,h'}(x_0) - K_{h'}\star (K_h\star \check f)(x_0) \big)^2 + \big(\hat f_{h'}(x_0) - K_{h'}\star \check f(x_0) \big)^2 \right. \\
&\left. \hspace{2.7in} + \big( K_{h'}\star (K_h\star \check f)(x_0) - K_{h'}\star \check f(x_0)  \big)^2 - \frac{V(x_0,h')}{3}\right\}_+ \\
&\le 3\left( B(h) + D_1 + D_2 \right),
\end{align*}
where
\begin{align*}
B(h) &= \umax{h'\in\Hn} \left( K_{h'}\star (K_h\star \check f)(x_0) - K_{h'}\star \check f(x_0)  \right)^2 \\
D_1 &= \umax{h'\in\Hn}\left\{ \big(\hat f_{h'}(x_0) - K_{h'}\star \check f(x_0) \big)^2 - \frac{V(x_0,h')}{6} \right\}_+ \\
D_2 &= \umax{h'\in\Hn}\left\{ \big( \hat f_{h,h'}(x_0) - K_{h'}\star (K_h\star \check f)(x_0) \big)^2 - \frac{V(x_0,h')}{6} \right\}_+ .
\end{align*}
Since 
\begin{align}
B(h) &= \umax{h'\in\Hn}\left( K_{h'}\star (K_h\star \check f)(x_0) - K_{h'}\star \check f(x_0)  \right)^2 = \umax{h'\in\Hn} \left(K_{h'}\star (K_h\star \check f - \check f)(x_0) \right)^2 
 \notag \\
&\le \norme{K}{1}^2 \usup{t\in \mathcal{V}_n(x_0)}\big|K_h\star\check f(t) - \check f(t) \big|^2, \notag
\end{align}
then we can rewrite \eqref{eq:bound-fhat} as
\begin{multline}\label{eq:bound-fhat-bis}
\left(\hat f(x_0) - \check f(x_0) \right)^2 \le 18D_1 + 18D_2 + 6V(x_0,h)  + 3\big(\hat f_h(x_0)  - K_h \star \check f(x_0) \big)^2  \\ + (18\norme{K}{1}^2 + 3)\usup{t\in \mathcal{V}_n(x_0)}\big|K_h\star\check f(t) - \check f(t) \big|^2.
\end{multline}
The last two terms of \eqref{eq:bound-fhat-bis} are controlled by \eqref{eq:upper-bound-variance-term} and \eqref{eq:control-bias-fcheck} of Proposition \ref{prop:control-bias-fcheck}. Hence it remains to deal with terms $D_1$ and $D_2$. 

For $D_1$, we recall that $K_{h}\star \check f(x_0) = \tilde\E\big[\hat f_h(x_0) \big]$ and 
\begin{align}
\tilde\E[D_1] &= \tilde\E \left[ \umax{h \in\Hn}\left\{ \left(\hat f_{h}(x_0) - K_{h}\star \check f(x_0) \right)^2 - \frac{V(x_0,h)}{6} \right\}_+ \right]  \notag \\
&\le \sum_{h\in \Hn} \tilde\E \left[ \left\{ \left(\hat f_{h}(x_0) - \tilde\E\big[\hat f_h(x_0)  \right)^2 - \frac{V(x_0,h)}{6} \right\}_+ \right] \notag  \\
&\le  \sum_{h\in \Hn} \int_0^{+\infty} \tilde\Prob \left( \left\{ \big(\hat f_{h}(x_0) - \tilde\E\big[\hat f_h(x_0)  \big)^2 - \frac{V(x_0,h)}{6} \right\}_+ > u \right)du \notag  \\
&\le \sum_{h\in \Hn}  \int_0^{+\infty} \tilde\Prob \left( \big|\hat f_{h}(x_0) - \tilde\E\big[\hat f_h(x_0)  \big| > \sqrt{\frac{V(x_0,h)}{6} + u}\,\right)du. \label{eq:bound-prob-D1}
\end{align}

Now let us introduce the sequence of \textit{i.i.d.} random variables $Z_1,\ldots, Z_n$ where we set 
\[
Z_i = w(\tilde\theta_n, \hat g(X_i))K_h(x_0 - X_i).
\]
Then we have 
\[
\hat f_{h}(x_0) - \tilde\E\big[\hat f_h(x_0)\big]  = \frac{1}{n}\sumi \big (Z_i - \tilde\E[Z_i] \big ).
\]
Moreover, we have by \eqref{eq:upperbound-w-theta-g-bis} and recall that we are on $\Omega_\rho = \left\{ \rho^{-1}\gamma \le \hat\gamma \le \rho\gamma \right\}$,
\[
|Z_i| = |w(\tilde\theta_n, \hat g(X_i))K_h(x_0 - X_i)| \le \frac{2(\hat\gamma + 1)\normsupglobal{K}}{h\delta\hat\gamma} \le \frac{2(\rho\gamma + 1)\normsupglobal{K}}{h\delta\hat\gamma} =: b,
\]
and
\[
\tilde\E\big[ Z_1^2\big] = \tilde\E\left[ w(\tilde\theta_n, \hat g(X_i))^2K_h^2(x_0 - X_i)\right] \le \frac{4(\rho\gamma + 1)^2\norme{K}{2}^2\normsup{g}{\mathcal{V}_n(x_0)}}{h\delta^2 \hat\gamma^2}  =: v.
\]
Applying the Bernstein inequality (cf. Lemma 2 of Comte and Lacour \cite{ComteLacour}), we have for any $u>0$,
\begin{align*}
&\tilde\Prob \left( \big|\hat f_{h}(x_0) - \tilde\E\big[\hat f_h(x_0)  \big| > \sqrt{\frac{V(x_0,h)}{6} + u}\,\right) = \tilde\Prob \left( \Big| \frac{1}{n}\sumi \big (Z_i - \tilde\E[Z_i] \big )  \Big| > \sqrt{\frac{V(x_0,h)}{6} + u}\,\right) \\
&\hspace{0.75in} \le 2\max\left\{\exp\left(- \frac{n}{4v}\left(\frac{V(x_0,h)}{6} + u\right) \right), \exp\left( - \frac{n}{4b}\sqrt{\frac{V(x_0,h)}{6} + u} \right) \right\}\\
&\hspace{0.75in} \le 2\max\left\{\exp\left(- \frac{n}{24v}V(x_0,h)\right)\exp\left( - \frac{nu}{4v}\right), \exp\left(-  \frac{n}{8b}\sqrt{\frac{V(x_0,h)}{6}} \right) \exp\left( - \frac{n\sqrt{u}}{8b}\right)\right\}
\end{align*}
On the other hand, by the definition of $V(x_0,h)$ we have
\begin{align*}
\frac{n}{24v}V(x_0,h) & =  \frac{nh\hat\gamma^2\delta^2}{96(\rho\gamma + 1)^2 \norme{K}{2}^2\normsup{g}{\mathcal{V}_n(x_0)}} \times  \frac{\kappa \norme{K}{1}^2\norme{K}{2}^2\normsup{g}{\mathcal{V}_n(x_0)}}{\hat\gamma^2 nh}\log(n) \\ 
& =  \frac{\kappa \delta^2 \norme{K}{1}^2}{96(\rho\gamma+1)^2}\log(n)  \ge  \dfrac{\kappa \delta^2}{96(\rho\gamma+1)^2}\log(n). 
\end{align*}
If we choose $\kappa$ such that  $\dfrac{\kappa \delta^2}{96(\rho\gamma+1)^2} \ge 2$, we get 
\[
\frac{n}{24v}V(x_0,h) \ge 2\log(n). 
\]
Moreover, using the assumption that $\hat\gamma nh \ge \log^3(n)$ and that $\hat\gamma\leq\rho\gamma$ on $\Omega_{\rho}$, we have
\begin{align*}
\frac{n}{8b}\sqrt{\frac{V(x_0,h)}{6}} &=   \frac{nh\hat\gamma \delta}{16\sqrt{6}(\rho\gamma+1)\normsupglobal{K}}  \times \frac{\norme{K}{1}\norme{K}{2}\sqrt{\kappa \normsup{g}{\mathcal{V}_n(x_0)}\log(n)}}{\hat\gamma \sqrt{nh}} \\
&=  \frac{\delta \norme{K}{1}\norme{K}{2}\normsup{g}{\mathcal{V}_n(x_0)}^{1/2} }{16\sqrt{6}(\rho\gamma+1)\normsupglobal{K}}\sqrt{\kappa nh \log(n)}  \\
&\ge  \frac{\delta \norme{K}{1}\norme{K}{2}}{16\sqrt{6}(\rho\gamma + 1)\rho^{1/2}\gamma^{1/2} \normsupglobal{K}}\sqrt{\kappa}\log^2(n) \ge 2\log(n),
\end{align*}
if 
\[
\frac{\delta \norme{K}{1}\norme{K}{2}}{16\sqrt{6}(\rho\gamma + 1)\rho^{1/2}\gamma^{1/2} \normsupglobal{K}}\sqrt{\kappa} \log(n) \ge 2 
\]
which automatically holds for well-chosen value of $\kappa$, and $n$ large enough. Then we have by using the conditions $\rho^{-1}\gamma \le \hat\gamma$ and $h\ge 1/n$, 
\begin{align*}
\tilde\E[D_1] &\le \sum_{h\in\Hn} \int_0^{+\infty} 2n^{-2} \max \left\{\exp\left( - \frac{nu}{4v}\right), \exp\left( - \frac{n\sqrt{u}}{8b}\right)\right\}du \\
&\le 2n^{-2}\sum_{h\in\Hn} \int_0^{+\infty} \max \left\{\exp\left( -  nh\frac{\delta^2 \hat\gamma^2}{16 (\rho\gamma+1)^2 \norme{K}{2}^2\normsup{g}{\mathcal{V}_n(x_0)}}u  \right), \exp\left( - nh \frac{\delta\hat\gamma}{16(\rho\gamma + 1) \normsupglobal{K}}\sqrt{u}  \right)\right\}du \\
&\le 2n^{-2}\sum_{h\in\Hn} \int_0^{+\infty} \max \left\{\exp\left( -  nh\frac{\delta^2\gamma^2}{16 (\rho\gamma+1)^2 \rho^2 \norme{K}{2}^2\normsup{g}{\mathcal{V}_n(x_0)}}u  \right), \exp\left( - nh \frac{\delta\gamma}{16(\rho\gamma + 1) \rho \normsupglobal{K}}\sqrt{u}  \right)\right\}du \\
&\le 2n^{-2}\sum_{h\in\Hn} \int_0^{+\infty} \max \left\{e^{- \pi_1 u}, e^{- \pi_2 \sqrt{u}}\right\}du \le 2n^{-2} \sum_{h\in\Hn} \max\left\{\frac{1}{\pi_1}, \frac{2}{\pi_2^2} \right\}.
\end{align*}
with $\pi_1:= \dfrac{\delta^2\gamma^2}{16 (\rho\gamma+1)^2 \rho^2 \norme{K}{2}^2\normsup{g}{\mathcal{V}_n(x_0)}}$ and $\pi_2:=  \dfrac{\delta\gamma}{16(\rho\gamma + 1) \rho \normsupglobal{K}}$. \smallskip

Since $\card(\Hn) \le n$, we finally obtain 
\begin{equation}\label{eq:D1-upperbound}
\tilde\E[D_1] \le C_5\delta^{-2}\gamma^{-2} n^{-1},
\end{equation}
where $C_5$ is a positive constant depending on $\normsup{g}{\mathcal{V}_n(x_0)}$, $\normsupglobal{K}$, $\norme{K}{2}$ and $\rho$. 

Similarly, we introduce $U_i = w(\tilde\theta_n,\hat g(X_i))K_{h'}\star K_h(x_0 - X_i)$ for $i=1,\ldots, n$. Then,
\[
\hat f_{h,h'}(x_0) - K_{h'}\star (K_h\star \check f)(x_0) = \hat f_{h,h'}(x_0) - \tilde\E\big[ \hat f_{h,h'}(x_0) \big] = \frac{1}{n} \sumi \big(U_i - \tilde\E[U_i]\big),
\]
and
\[
|U_i| \le \frac{4\norme{K}{1} \normsupglobal{K}}{h'\delta\hat\gamma} =: \bar b,\quad \text{ and } \quad \tilde\E\big[ U_1^2\big]  \le \frac{16 \norme{K}{1}^2\norme{K}{2}^2\normsup{g}{\mathcal{V}_n(x_0)}}{h'\delta^2\hat\gamma^2} =: \bar v.
\]
Following the same lines as for obtaining \eqref{eq:D1-upperbound}, we get by using Bernstein inequality
\begin{equation}\label{eq:D2:upperbound}
\tilde\E[D_2] \le C_6\delta^{-2}\gamma^{-2} n^{-1},
\end{equation}
with $C_6$ a positive constant depends on  $\normsup{g}{\mathcal{V}_n(x_0)}$, $\normsupglobal{K}$, $\norme{K}{1}$, $\norme{K}{2}$ and $\rho$. 

Finally, combining \eqref{eq:bound-fhat-bis},  \eqref{eq:D1-upperbound}, \eqref{eq:D2:upperbound} and successively applying Lemma \ref{lem:control-Omega-rho} and Lemma \ref{lem:control-risk-on-OmegaC} allow us to conclude the result stated in Theorem \ref{th:oracle-inequality}.

\subsection{Proof of Lemma \ref{lem:g_sym}}\label{sec:proof_lem_g_sym}
First, we prove that  $g^{sym}$ is the density of $Y_i$.  To this aim, let $\varphi$ be  a measurable bounded function defined on $\R$. We compute\begin{eqnarray*}\mathbb E[\varphi(Y_i)]&=&\mathbb E[\mathbb E[\varphi(X_i)|\varepsilon_i]\id_{\{\epsilon_1 = 1\}}]+\mathbb E[\mathbb E[\varphi(2-X_i)|\varepsilon_i]\id_{\{\epsilon_1 = -1\}}],\\
&=&\frac{1}{2}\left(\mathbb E[\varphi(X_i)] +\mathbb E[\varphi(2-X_i)]\right),\\
&=&\frac{1}{2}\left(\int_0^1\varphi(x)g(x)dx +\int_0^1\varphi(2-x)g(x)dx\right),\\
&=&\frac{1}{2}\left(\int_0^1\varphi(x)g(x)dx +\int_1^2\varphi(x)g(2-x)dx\right),\\
&=&\int_0^2\varphi(x)g^{sym}(x)dx.\end{eqnarray*}
Since the equality holds for any test function $\varphi$, we obtain the first assumption of the lemma.

We prove now \eqref{eq:theta-upper-bound}. Under the identifiability condition, we have $\theta = g(x)$ for all $x\in [1-\delta, 1]$, and thus  $\theta=2g^{sym}(x)$ for $x\in [1-\delta,1+\delta]$. Hence we have

\begin{align*}
|\hat\theta_{n,b} - \theta| &=\left| \frac{1}{\delta}\int_{1-\delta}^{1+\delta} \hat{g}^{sym}_b(x)dx - \frac{1}{\delta}\int_{1-\delta}^{1+\delta} g^{sym}(x) dx \right| = \left| \frac{1}{\delta}\int_{1-\delta}^{1+\delta} \hat{g}^{sym}_b(x) -g^{sym}(x)dx \right| \\
&\leq \frac{1}{\delta} \int_{1-\delta}^{1+\delta} \left| \hat{g}^{sym}_b(x) - g^{sym}(x) \right| dx  \\
&\leq \frac{1}{\delta} \int_{1-\delta}^{1+\delta} \normsup{\hat{g}^{sym}_b - g^{sym}}{[1-\delta,1+\delta]} dx = 2\normsup{\hat{g}^{sym}_b - g^{sym}}{[1-\delta,1+\delta]},
\end{align*}


which proves \eqref{eq:theta-upper-bound}. Then, thanks to the Markov Inequality
\begin{eqnarray*}
\Prob\left(\tilde\theta_{n,b}\neq \hat\theta_{n,b}\right)&=&\Prob\left(\hat\theta_{n ,b}\notin \left[\frac{\delta}{2}, 1-\frac{\delta}{2}\right] \right) \\
&\leq& \Prob\left( |\hat{\theta}_{n,b} - \theta| > \frac{\delta}{2}\right) \leq \frac{4}{\delta^2}\E\left[ |\hat\theta_{n,b} - \theta|^2 \right],
\end{eqnarray*}
which is \eqref{eq:markov}. Finally, 
\begin{align*}
\E \left[|\tilde{\theta}_{n,b} - \theta|^2 \right] &=\E \left[|\tilde{\theta}_{n,b} - \theta|^2\left(\mathbf 1_{\big\{\hat\theta_{n,b} = \tilde\theta_{n,b}\big\}}+\mathbf 1_{\big\{\hat\theta_{n,b} \neq \tilde\theta_{n,b}\big\}}\right) \right]\\
&\le \E \left[|\hat{\theta}_{n,b} - \theta|^2 \mathbf 1_{\big\{\hat\theta_{n,b} \in [\delta/2,1-\delta/2]\big\}}\right] + \left(|\tilde{\theta}_{n,b}| + |\theta|\right)^2\mathbb P\left(\hat\theta_{n,b} \ne \tilde\theta_{n,b}\right)
\\
&\le \E \left[|\hat{\theta}_{n,b} - \theta|^2 \mathbf 1_{\big\{\hat\theta_{n,b} \in [\delta/2,1-\delta/2]\big\}}\right] + 4\mathbb P\left(\hat\theta_{n,b} \ne \tilde\theta_{n,b}\right)\\
&\le (1+4\times \frac{4}{\delta^2})\E \left[|\hat{\theta}_{n,b} - \theta|^2 \right],\\
&\leq (1+4\times \frac{4}{\delta^2}) \times 2\E \left[\normsup{\hat{g}_b^{sym} - g^{sym}}{[1-\delta,1+\delta]}^2\right]. 
\end{align*}
thanks to  \eqref{eq:markov} and then \eqref{eq:theta-upper-bound}. This concludes the proof of Lemma \ref{lem:g_sym}.

\subsection{Proof of Corollary \ref{cor:rate-cv-fhat} }

Since Assumptions {\bf (A6)} and {\bf (A7)} are fulfilled. According to Proposition 1.2 of Tsybakov \cite{Tsybakov04}, we get for all $x_0\in [0,1]$
\[
|K_h\star f(x_0) - f(x_0)| \le C_7 \mathcal{L} h^{\beta},
\]
where $C$ a constant depending on $K$ and $\mathcal{L}$. 
We obtain
\begin{equation}\label{eq:rate-cv-bias-variance}
\min_{h\in\mathcal H_n}\left\{\normsup{K_h\star f - f}{\mathcal{V}_n(x_0)}^2 + \frac{\log(n)}{\delta^2\gamma^2 nh}\right\} \le \min_{h\in\mathcal H_n}\left\{C_7\mathcal L h^\beta+\frac{\log(n)}{\delta^2\gamma^2 nh}\right\}. 
\end{equation}
Taking
\[
h^* = \frac1{k^*} \text{ with } k^*=\left\lfloor \left(\frac{n}{\log(n)}\right)^{1/(2\beta+1)}\right\rfloor,
\]
there exists $n(\beta,\gamma,\rho)$ such that, for all $n\geq n(\beta,\gamma,\rho)$,
\[
\frac{\gamma}\rho\frac{n}{\log^3(n)}\geq\left\lfloor\left(\frac{n}{\log(n)}\right)^{1/(2\beta+1)}\right\rfloor\geq \log(n)=\alpha_n. 
\]
Implying that for all $n\geq n(\beta,\gamma,\rho)$, 
\[
\Omega_\rho\subseteq\{\hat\gamma\geq\gamma/\rho\}\subseteq\{h^*\in\mathcal H_n\}. 
\]

Finally, since we also have \eqref{eq:rate-cv-g} and \eqref{eq:rate-cv-theta}, gathering \eqref{eq:oracle-inequality} and \eqref{eq:rate-cv-bias-variance}. Since Assumption {\bf (A5)} is verified by construction of $\mathcal H_n$, using again Lemma~\ref{lem:control-risk-on-OmegaC}, we obtain, for all $n$,
\[
 \E\left[\big(  \hat f(x_0) - f(x_0)\big)^2\right] \le C_8 \left( \frac{\log n}{n}\right)^{\frac{2\beta}{2\beta + 1}},
\]
where $C_8$ is a constant depending on $K$, $\normsup{f}{\mathcal{V}_n(x_0)}$, $g$, $\delta$, $\gamma$, $\rho$, $\mathcal{L}$ and $\beta$. 

\subsection{Proofs of technical intermediate results}\label{sec:intermediate_proof}

\subsubsection{Proof of Proposition \ref{prop:control-bias-fcheck}}\label{sec:proof-prop-control-bias}
Let us introduce the function 
\begin{equation}\label{eq:f-tilde}
\tilde f(x) := w(\tilde\theta_n, g(x))g(x) = \frac{1}{1-\tilde\theta_n}\left(1 - \frac{\tilde\theta_n}{g(x)} \right)g(x).
\end{equation}
Then we have for $x_0\in [0,1]$
\[
\big( \check f(x_0) - f(x_0) \big)^2 \le 2\left( \big( \check f(x_0) - \tilde f(x_0) \big)^2 + \big(\tilde f(x_0) - f(x_0) \big)^2 \right).
\]
For the first term, on $\Omega_\rho = \left\{ \rho^{-1}\gamma \le \hat\gamma \le \rho\gamma \right\}$ we have, by using {\bf (A4)},
\begin{align}
\big( \check f(x_0) - \tilde f(x_0) \big)^2 &=  \left(w(\tilde\theta_n, \hat g(x_0))g(x_0)  - w(\tilde\theta_n, g(x_0))g(x_0) \right)^2 \notag\\ &= \left(\frac{1}{1 -\tilde\theta_n}\Bigg(1 - \frac{\tilde\theta_n}{\hat g(x_0)} \Bigg) - \frac{1}{1 -\tilde\theta_n}\Bigg(1 - \frac{\tilde\theta_n}{g(x_0)} \Bigg)\right)^2|g(x_0)|^2 \notag  \\
&=\frac{\tilde\theta_n^2}{(1 - \tilde\theta_n)^2}   \left(\frac{1}{\hat g(x_0)}  - \frac{1}{g(x_0)}\right)^2|g(x_0)|^2\notag \\
&\le \frac{4}{\delta^2}\left(\frac{\hat g(x_0) - g(x_0)}{\hat g(x_0) g(x_0)}\right)^2 |g(x_0)|^2 \notag \\
&\le 4\rho^2\delta^{-2} \gamma^{-2} \normsup{\hat g - g}{\mathcal{V}_n(x_0)}^2.  \label{eq:upperbound-fcheck-ftilde}
\end{align}

Moreover, thanks to {\bf(A1)},
\begin{align}
\big( \tilde f(x_0) - f(x_0) \big)^2 &= \left(w(\tilde{\theta}_n, g(x_0))g(x_0)  - w(\theta, g(x_0))g(x_0) \right)^2 \notag \\
&= \left(\frac{1}{1 - \tilde\theta_n}\Bigg(1 - \frac{\tilde{\theta}_n}{g(x_0)} \Bigg)g(x_0)	- \frac{1}{1-\theta}\Bigg(1 - \frac{\theta}{g(x_0)} \Bigg)g(x_0) \right)^2 \notag \\
&= \left( \frac{1}{1 - \tilde{\theta}_n} - \frac{1}{1-\theta}	+ \Bigg(\frac{\theta}{1-\theta}	-  \frac{\tilde{\theta}_n}{1 - \tilde\theta_n} \Bigg) \frac{1}{g(x_0)} \right)^2 |g(x_0)|^2\notag \\
&= \frac {|g(x_0)|^2}{(1-\theta)^2(1-\tilde{\theta}_n)^2} \left( \tilde{\theta}_n - \theta + \frac{\theta - \tilde{\theta}_n}{g(x_0)} \right)^2 \notag \\
&\le \frac{4\normsup{g}{\mathcal{V}_n(x_0)}^2}{\delta^4} \left( \tilde{\theta}_n - \theta + \frac{\theta - \tilde{\theta}_n}{g(x_0)}\right)^2 \notag \\
&\le 16\normsup{g}{\mathcal{V}_n(x_0)}^2\delta^{-6}\big|\tilde{\theta}_n - \theta\big|^2. \label{eq:upperbound-ftilde-f}
\end{align}
Thus we obtain by gathering \eqref{eq:upperbound-fcheck-ftilde} and \eqref{eq:upperbound-ftilde-f},
\[
\big( \check f(x_0) - f(x_0) \big)^2 \le 8\rho^2\delta^{-2}\gamma^{-2} \normsup{\hat g - g}{\mathcal{V}_n(x_0)}^2 + 32\normsup{g}{\mathcal{V}_n(x_0)}^2\delta^{-6}\big|\tilde{\theta}_n - \theta\big|^2.
\]
Next, the term $\big( K_h\star \check f(x_0) - \check f(x_0)\big)^2$ can be treated by studying the following decomposition
\begin{align*}
\big( K_h\star \check f(x_0) - \check f(x_0)\big)^2 &\le 3\bigg(\big( K_h\star \check f(x_0) - K_h\star \tilde f(x_0) \big)^2 
+ \big( K_h\star \tilde f(x_0) - K_h\star f(x_0) \big)^2  \\
&\hspace{2.7in} + \big( K_h\star  f(x_0) - \check f(x_0) \big)^2\bigg) \\&=: 3\big(A_1 + A_2 + A_3).
\end{align*}
For term $A_1$, we have by using \eqref{eq:upperbound-fcheck-ftilde}
\begin{align}
A_1 = \big( K_h\star (\check f - \tilde f)(x_0) \big)^2 &= \left(\int K_h(x_0 - u)(\check f(u) - \tilde f(u))du \right)^2 \notag \\
&\le \left(\int |K_h(x_0 - u)||\check f(u) - \tilde f(u)|du \right)^2 \notag \\
&\le 4\rho^2\delta^{-2}\gamma^{-2}  \normsup{\hat g - g}{\mathcal{V}_n(x_0)}^2\left( \int |K_h(x_0 - u)|du \right)^2\notag \\
&\le 4\rho^2\delta^{-2}\gamma^{-2} \norme{K}{1}^2 \normsup{\hat g - g}{\mathcal{V}_n(x_0)}^2 . \notag
\end{align}
By using \eqref{eq:upperbound-ftilde-f} and following the same lines as for $A_1$, we obtain
\[
A_2 = \big( K_h\star (\tilde f- f)(x_0) \big)^2 \le 16\normsup{g}{\mathcal{V}_n(x_0)}^2\delta^{-6}\norme{K}{1}^2 \big|\tilde{\theta}_n - \theta\big|^2.
\]
For $A_3$, using the upper bound obtained as above for $( \check f(x_0) - f(x_0) )^2$, we have
\begin{align*}
A_3 &\le 2\big( K_h\star f(x_0) - f(x_0) \big)^2 + 2\big( f(x_0) - \check f(x_0)  \big)^2 \\
&\le 2\normsup{K_h\star f - f}{\mathcal{V}_n(x_0)}^2 + 16\rho^2\delta^{-2}\gamma^{-2} \normsup{\hat g - g}{\mathcal{V}_n(x_0)}^2 + 64\normsup{g}{\mathcal{V}_n(x_0)}^2\delta^{-6}\big|\tilde{\theta}_n - \theta\big|^2.
\end{align*}
Finally, combining all the terms $A_1$, $A_2$ and $A_3$, we obtain \eqref{eq:control-bias-fcheck}. This ends the proof of Proposition \ref{prop:control-bias-fcheck}.

\subsubsection{Proof of Lemma \ref{lem:control-Omega-rho}}
Lemma \ref{lem:control-Omega-rho} is a consequence of \eqref{eq:condtion-on-ghat}. Indeed, if condition {\bf (A3)}  is satisfied,  we have for all $t\in \mathcal{V}_n(x_0)$, $|\hat g(t) - g(t)|\le \nu|\hat g(t)|$ with probability $1- C_{g,\nu} \exp\big(-(\log n)^{3/2}\big)$. 

This implies,
\[
(1+\nu)^{-1}|g(t)| \le |\hat g(t)| \le (1-\nu)^{-1}|g(t)|.
\]
Since $\gamma = \underset{t\in \mathcal{V}_n(x_0)}{\inf} |g(t)|$ and $\hat\gamma = \underset{t\in \mathcal{V}_n(x_0)}{\inf} |\hat g(t)|$, by using \eqref{eq:condtion-on-ghat} and taking $\nu = \rho - 1$, $\nu = 1-\rho^{-1}$, we obtain with probability $1- C_{g,\nu} \exp\big(-(\log n)^{3/2}\big)$,  $(1+\nu)^{-1}\gamma \le \hat\gamma \le (1-\nu)^{-1}\gamma$. This completes the proof of Lemma \ref{lem:control-Omega-rho}. 

\subsubsection{Proof of Lemma \ref{lem:control-risk-on-OmegaC}}
We have for any $x_0\in [0,1]$,
\begin{align*}
\E\left[ \big( \hat f_h(x_0) - f(x_0) \big)^2 \id_{\Omega_\rho^c}\right] \le 2\E\big[|\hat f_h(x_0)|^2 \id_{\Omega_\rho^c}\big] + 2\normsup{f}{\mathcal{V}_n(x_0)}^2\Prob(\Omega_\rho^c).
\end{align*}
\begin{align}
\E\big[|\hat f_h(x_0)|^2 \id_{\Omega_\rho^c}\big] &= \E\left[\Bigg|\frac{1}{nh}\sumi w(\tilde\theta_n,\hat g(X_i) )K\left(\frac{x_0 - X_i}{h} \right)\Bigg|^2 \id_{\Omega_\rho^c} \right] \notag \\
&\le \E\left[\left(\frac{1}{nh}\sumi \left|w(\tilde\theta_n,\hat g(X_i) )K\left(\frac{x_0 - X_i}{h}\right)\right| \right)^2 \id_{\Omega_\rho^c} \right] \notag \\
&\le  \frac{4}{\delta^2}\E\left[ \left(\frac{1}{nh} \sumi \left|K\left(\frac{x_0 - X_i}{h}\right)\right| \right)^2 \left(1+\frac{1}{\hat \gamma}\right)^2 \id_{\Omega_\rho^c} \right] \text{ (using \eqref{eq:upperbound-w-theta-g-bis})}  \notag \\
&\le  \frac{4}{\delta^2}\E\left[ \left(\frac{1}{nh} \sumi \left|K\left(\frac{x_0 - X_i}{h}\right)\right| \right)^2 \right] \E\left[ \left(1+\frac{1}{\hat \gamma}\right)^2 \id_{\Omega_\rho^c} \right] \text{ (by independence)}  \notag \\
&\le  \frac{4\normsupglobal{K}^2}{\delta^2 h^2} \E\left[ \left(1+\frac{1}{\hat \gamma}\right)^2 \id_{\Omega_\rho^c} \right]  \notag \\
&\le \frac{4\normsupglobal{K}^2}{\delta^2} n^2 \left(1+\frac{1}{(\log n)^3}\right)^2 \Prob(\Omega_\rho^c) \text{ (using Assumption {\bf (A5)}).}   \notag
\end{align}
Finally, we apply Lemma \ref{lem:control-Omega-rho} to establish the following bound
\begin{align*}
\E\left[ \big( \hat f_h(x_0) - f(x_0) \big)^2 \id_{\Omega_\rho^c}\right] &\le C_{g,\rho}\left( \frac{8\normsupglobal{K}^2}{\delta^2} \frac{n^2}{(\log n)^6} + 2\normsup{f}{\mathcal{V}_n(x_0)}^2 \right)\exp\left\{ - (\log n)^{3/2} \right\}  \\
&\le \frac{C}{n^2},
\end{align*}
where $C$ depends on $\delta$, $\normsup{f}{\mathcal{V}_n(x_0)}$, $\normsupglobal{K}$, $g$ and $\rho$, which ends the proof of Lemma \ref{lem:control-risk-on-OmegaC}. 

\subsection*{Acknowledgment}
We are very grateful to Catherine Matias for interesting discussions on mixture models. The research of the authors is partly supported by the french Agence Nationale de la Recherche (ANR-18-CE40-0014 projet SMILES) and by the french R\'egion Normandie (projet RIN AStERiCs 17B01101GR). Finally, we gratefully acknowledge the referees for carefully reading the manuscript and for numerous suggestions that improved the paper.

\bibliographystyle{plain}	
\bibliography{references} 

\def\cprime{$'$}
\begin{thebibliography}{10}

\bibitem{Benjamini1995}
Yoav Benjamini and Yosef Hochberg.
\newblock Controlling the false discovery rate: a practical and powerful
  approach to multiple testing.
\newblock {\em Journal of the Royal statistical society: series B
  (Methodological)}, 57(1):289--300, 1995.

\bibitem{BLR2013}
Karine Bertin, Claire Lacour, and Vincent Rivoirard.
\newblock Adaptive pointwise estimation of conditional density function.
\newblock {\em preprint arXiv:1312.7402}, 2013.

\bibitem{BLR2016}
Karine Bertin, Claire Lacour, and Vincent Rivoirard.
\newblock Adaptive pointwise estimation of conditional density function.
\newblock {\em Ann. Inst. H. Poincaré Probab. Statist.}, 52(2):939--980, 05
  2016.

\bibitem{Butucea2000}
Cristina Butucea.
\newblock Two adaptive rates of convergence in pointwise density estimation.
\newblock {\em Math. Methods Statist.}, 9(1):39--64, 2000.

\bibitem{celisse2010}
Alain Celisse and St{\'e}phane Robin.
\newblock A cross-validation based estimation of the proportion of true null
  hypotheses.
\newblock {\em Journal of Statistical Planning and Inference},
  140(11):3132--3147, 2010.

\bibitem{Chagny13}
Ga{\"e}lle Chagny.
\newblock Penalization versus goldenshluger- lepski strategies in warped bases
  regression.
\newblock {\em ESAIM: Probability and Statistics}, 17:328--358, 2013.

\bibitem{Chichignoud2017}
Micha{\"e}l Chichignoud, Van~Ha Hoang, Thanh Mai~Pham Ngoc, Vincent Rivoirard,
  et~al.
\newblock Adaptive wavelet multivariate regression with errors in variables.
\newblock {\em Electronic journal of statistics}, 11(1):682--724, 2017.

\bibitem{comte_estimation_2015}
Fabienne Comte.
\newblock {\em Estimation non-param{\'e}trique}.
\newblock Spartacus-IDH, 2015.

\bibitem{comte_adaptive_2011}
Fabienne Comte, Stéphane Ga\"{\i}ffas, and Agathe Guilloux.
\newblock Adaptive estimation of the conditional intensity of marker-dependent
  counting processes.
\newblock {\em Ann. Inst. Henri Poincar\'{e} Probab. Stat.}, 47(4):1171--1196,
  2011.

\bibitem{ComteSamson13}
Fabienne Comte, Valentine Genon-Catalot, and Adeline Samson.
\newblock Nonparametric estimation for stochastic differential equations with
  random effects.
\newblock {\em Stochastic Processes and their Applications}, 123(7):2522--2551,
  2013.

\bibitem{ComteLacour}
Fabienne Comte and Claire Lacour.
\newblock Anisotropic adaptive kernel deconvolution.
\newblock {\em Ann. Inst. Henri Poincar\'e Probab. Stat.}, 49(2):569--609,
  2013.

\bibitem{ComteRebafka2016}
Fabienne Comte and Tabea Rebafka.
\newblock Nonparametric weighted estimators for biased data.
\newblock {\em Journal of Statistical Planning and Inference}, 174:104--128,
  2016.

\bibitem{DHRR}
Marie Doumic, Marc Hoffmann, Patricia Reynaud-Bouret, and Vincent Rivoirard.
\newblock Nonparametric estimation of the division rate of a size-structured
  population.
\newblock {\em SIAM Journal on Numerical Analysis}, 50(2):925--950, 2012.

\bibitem{efron2001empirical}
Bradley Efron, Robert Tibshirani, John~D Storey, and Virginia Tusher.
\newblock Empirical bayes analysis of a microarray experiment.
\newblock {\em Journal of the American statistical association},
  96(456):1151--1160, 2001.

\bibitem{genovese2002operating}
Christopher Genovese and Larry Wasserman.
\newblock Operating characteristics and extensions of the false discovery rate
  procedure.
\newblock {\em Journal of the Royal Statistical Society: Series B (Statistical
  Methodology)}, 64(3):499--517, 2002.

\bibitem{gine_exponential_2009}
Evarist Gin\'{e} and Richard Nickl.
\newblock An exponential inequality for the distribution function of the kernel
  density estimator, with applications to adaptive estimation.
\newblock {\em Probab. Theory Related Fields}, 143(3-4):569--596, 2009.

\bibitem{GL11}
Alexander Goldenshluger and Oleg Lepski.
\newblock Bandwidth selection in kernel density estimation: orcale inequalities
  and adaptive minimax optimality.
\newblock {\em The Annals of Statistics}, 39(3):1608--1632, 2011.

\bibitem{huber1965robust}
Peter~J Huber.
\newblock A robust version of the probability ratio test.
\newblock {\em The Annals of Mathematical Statistics}, pages 1753--1758, 1965.

\bibitem{ibragimov_estimate_1980}
Ildar~A. Ibragimov and Rafael~Z. Has\cprime~minski\u{\i}.
\newblock An estimate of the density of a distribution.
\newblock {\em Zap. Nauchn. Sem. Leningrad. Otdel. Mat. Inst. Steklov. (LOMI)},
  98:61--85, 161--162, 166, 1980.
\newblock Studies in mathematical statistics, IV.

\bibitem{langaas2005}
Mette Langaas, Bo~Henry Lindqvist, and Egil Ferkingstad.
\newblock Estimating the proportion of true null hypotheses, with application
  to dna microarray data.
\newblock {\em Journal of the Royal Statistical Society: Series B (Statistical
  Methodology)}, 67(4):555--572, 2005.

\bibitem{Lepski2013_sup-norm-loss}
Oleg Lepski.
\newblock Multivariate density estimation under sup-norm loss: oracle approach,
  adaptation and independence structure.
\newblock {\em The Annals of Statistics}, 41(2):1005--1034, 2013.

\bibitem{liu2017density}
Haoyang Liu and Chao Gao.
\newblock Density estimation with contaminated data: Minimax rates and theory
  of adaptation.
\newblock {\em arXiv preprint arXiv:1712.07801}, 2017.

\bibitem{nguyen_matias_2014}
Van~Hanh Nguyen and Catherine Matias.
\newblock Nonparametric estimation of the density of the alternative hypothesis
  in a multiple testing setup. application to local false discovery rate
  estimation.
\newblock {\em ESAIM: Probability and Statistics}, 18:584–612, 2014.

\bibitem{Ng-Matias-SJS}
Van~Hanh Nguyen and Catherine Matias.
\newblock On efficient estimators of the proportion of true null hypotheses in
  a multiple testing setup.
\newblock {\em Scandinavian Journal of Statistics}, 41(4):1167--1194, 2014.

\bibitem{R-B14}
Patricia Reynaud-Bouret, Vincent Rivoirard, Franck Grammont, and Christine
  Tuleau-Malot.
\newblock Goodness-of-fit tests and nonparametric adaptive estimation for spike
  train analysis.
\newblock {\em The Journal of Mathematical Neuroscience}, 4(1):1, 2014.

\bibitem{Robin07}
St{\'{e}}phane Robin, Avner Bar{-}Hen, Jean{-}Jacques Daudin, and Laurent
  Pierre.
\newblock A semi-parametric approach for mixture models: Application to local
  false discovery rate estimation.
\newblock {\em Computational Statistics {\&} Data Analysis}, 51(12):5483--5493,
  2007.

\bibitem{Schuster85}
Eugene~F Schuster.
\newblock Incorporating support constraints into nonparametric estimators of
  densities.
\newblock {\em Communications in Statistics-Theory and methods},
  14(5):1123--1136, 1985.

\bibitem{storey2002}
John~D Storey.
\newblock A direct approach to false discovery rates.
\newblock {\em Journal of the Royal Statistical Society: Series B (Statistical
  Methodology)}, 64(3):479--498, 2002.

\bibitem{strimmer2008unified}
Korbinian Strimmer.
\newblock A unified approach to false discovery rate estimation.
\newblock {\em BMC bioinformatics}, 9(1):303, 2008.

\bibitem{Tsybakov04}
Alexandre~B Tsybakov.
\newblock {\em Introduction to Nonparametric Estimation}.
\newblock Springer series in Statistics. Springer, New York, 2009.

\end{thebibliography}

\end{document}